\documentclass[final]{siamltex}
\usepackage{amsmath}
\usepackage{amsfonts}
\usepackage{amssymb}

\usepackage{graphicx}

\usepackage{url}
\usepackage{xcolor}
\usepackage{fullpage}
\usepackage{bm}
\usepackage{mathrsfs}

\newcommand{\om}{\omega}

\newcommand{\la}{\lambda}
\newcommand{\ep}{\varepsilon}

\newcommand{\cY}{\mathcal{Y}}
\newcommand{\cB}{\mathcal{B}}
\newcommand{\cA}{\mathcal{A}}
\newcommand{\cS}{\mathcal{S}}
\newcommand{\cE}{\mathcal{E}}
\newcommand{\cU}{\mathcal{U}}
\newcommand{\cH}{\mathcal{H}}
\newcommand{\cT}{\mathcal{T}}
\newcommand{\cI}{\mathcal{I}}
\newcommand{\cN}{\mathscr{N} \hspace{-0.02in}}
\newcommand{\cP}{\mathcal{P}}
\newcommand{\sL}{\mathscr{L}}

\newcommand{\D}{\mathscr{D}}
\newcommand{\sS}{\mathscr{S}}

\newcommand{\brho}{\bm{\rho}}
\newcommand{\bg}{{\bf g}}
\newcommand{\bG}{\mathbf{\mathcal G}}
\newcommand{\bd}{{\bf d}}
\newcommand{\bu}{{\bf u}}
\newcommand{\bw}{{\bf w}}
\newcommand{\bv}{{\bf v}}

\newcommand{\vy}{\vec{\mathbf{y}}}
\newcommand{\vx}{\vec{\mathbf{x}}}
\newcommand{\vz}{\vec{\mathbf{z}}}
\newcommand{\vbh}{\vec{\bf h}}

\newcommand{\bz}{{\bf z}}
\newcommand{\bx}{{\bf x}}
\newcommand{\vzeta}{\vec{\bm{\zeta}}}

\newcommand{\bxi}{{\bm{\xi}}}
\newcommand{\balpha}{{\bm{\alpha}}}

\newcommand{\real}{\mathbb{R}}

\newcommand{\complex}{\mathbb{C}}

\newcommand{\lin}{\left <}
\newcommand{\rin}{\right >}

\newtheorem{remark}[theorem]{Remark} 

\begin{document}

\title{Resolution analysis of imaging with $\ell_1$ optimization}
\author{Liliana Borcea \and Ilker Kocyigit \footnotemark[1]}
\renewcommand{\thefootnote}{\fnsymbol{footnote}}
\footnotetext[1]{Department of Mathematics, University of Michigan,
  Ann Arbor, MI 48109-1043. \\\hspace{0.3in}Email: borcea@umich.edu \&
  ilkerk@umich.edu} \maketitle \date{today}
\begin{abstract}
  We study array imaging of a sparse scene of point-like sources or
  scatterers in a homogeneous medium. For source imaging the sensors
  in the array are receivers that collect measurements of the wave
  field. For imaging scatterers the array probes the medium with waves
  and records the echoes.  In either case the image formation is
  stated as a sparsity promoting $\ell_1$ optimization problem, and
  the goal of the paper is to quantify the resolution. We consider
  both narrow-band and broad-band imaging, and a geometric setup with
  a small array. We take first the case of the unknowns lying on the
  imaging grid, and derive resolution limits that depend on the
  sparsity of the scene.  Then we consider the general case with the
  unknowns at arbitrary locations. The analysis is based on estimates
  of the cumulative mutual coherence and a related concept, which we
  call interaction coefficient. It complements recent results in
  compressed sensing by deriving deterministic resolution limits that
  account for worse case scenarios in terms of locations of the
  unknowns in the imaging region, and also by interpreting the results
  in some cases where uniqueness of the solution does not hold. We
  demonstrate the theoretical predictions with numerical simulations.
  
\end{abstract}

\section{{Introduction}}
\label{sect:intro} 
Array imaging is an inverse problem for the wave equation, where the
goal is to determine remote sources or scatterers from measurements of
the wave field at a collection of nearby sensors, called the
array. The problem has applications in medical imaging, nondestructive
evaluation of materials, oil prospecting, seismic imaging, radar
imaging, ocean acoustics and so on. There is extensive literature on
various imaging approaches such as reverse time migration and its high
frequency version called Kirchhoff migration
\cite{Biondi,bleistein2001mathematics,CLAERBOUT}, matched field imaging 
\cite{baggeroer1993overview}, Multiple Signal
Classification (MUSIC) \cite{schmidt1986multiple,gruber2004time}, the
linear sampling method \cite{cakoni2003mathematical}, and the
factorization method \cite{kirsch1998characterization}.
In this paper we consider array imaging using $\ell_1$ optimization,
which is appropriate for sparse scenes of unknown sources or
scatterers that have small support in the imaging region. 

Imaging with sparsity promoting optimization has received much
attention recently, specially in the context of compressed sensing
\cite{fannjiang2010compressed,fannjiang2010compressive,rudelson2008sparse},
where a random set of sensors collect data from a sparse scene.  Such
studies use the restricted isometry property of the sensing matrix
\cite{candes2005decoding} or its mutual coherence
\cite{bruckstein2009sparse} to derive probability bounds on the event
that the imaging scene is recovered exactly for noiseless data, or
with small error that scales linearly with the noise level. The array
does not play an essential role in these studies, aside from its
aperture bounding the random sample of locations of the sensors, and
for justifying the scaling that leads to models of wave propagation
like the paraxial one \cite{fannjiang2010compressed}.

A different approach proposed in \cite{chai2013robust,chai2014imaging}
images a sparse scattering scene using illuminations derived from the
singular value decomposition (SVD) of the response matrix measured by
probing sequentially the medium with pulses emitted by one sensor at a
time and recording the echoes. Iluminations derived from the SVD are
known to be useful in imaging
\cite{prada1996decomposition,borcea2007optimal,borcea2010universal,
  borcea2008edge,gruber2004time} and they may mitigate noise. The
setup in \cite{chai2013robust}, which is typical in array imaging,
lets the sensors be closely spaced so that sums over them can be
approximated by integrals over the array aperture.  We consider the
same continuous aperture setup here and study the resolution of the
images produced by $\ell_1$ optimization, also known as basis
pursuit and
$\ell_1-$penalty. We address two questions: (1) How should we chose the
discretization of the imaging region so that we can guarantee unique
recovery of the sparse scene, at least when the unknowns lie on the
grid? (2) If the imaging region is discretized on a finer grid, for
which uniqueness does not hold, are there cases where the solution of
the $\ell_1$ optimization is still useful?

By studying question (1) we complement the existing results with
deterministic resolution limits that account for worse case scenarios,
and guarantee unique recovery of the scene for a given sparsity
$s$. This is defined as the number of non-zero entries of the vector
of unknowns or, equivalently, the number of grid points in the support
of the sources/scatterers that we image. We consider a geometric setup
with a small array, where wave propagation can be modeled by the
paraxial approximation. We have a more general paraxial model than in
\cite{fannjiang2010compressed}, which takes into consideration
sources/scatterers at different ranges from the array. This turns out
to be important in narrow-band regimes. We also consider broad-band
regimes and show that the additional multi-frequency data improves the
resolution.

It is typical in imaging with sparsity promoting optimization to
assume that the unknown sources or scatterers lie on the discretization
grid, meaning that they can be modeled by a sparse complex vector
$\brho \in \mathbb{C}^N$, where $N$ is the number of grid points.  If
the unknowns lie off-grid the results deteriorate. We refer to
\cite{herman2010general,fannjiang2013compressive} for a perturbation
analysis of compressed sensing with small off-grid displacements.
General tight error bounds can be found in
\cite{chi2011sensitivity}. They may be quite large and increase with
$N$.  Thus, there is a trade-off in imaging with $\ell_1$
optimization: on one hand we need a coarse enough discretization of
the imaging region to ensure unique recovery of the solution, and on
the other hand finer discretization to minimize modeling errors due to
off-grid placement of the unknowns. This trade-off is particularly
relevant in the narrow-band paraxial regime, where the resolution
limits may grow significantly with the sparsity $s$ of the scene.

At question (2) we consider fine discretizations of the imaging
region, to mitigate the modeling error. The problem is then how to
interpret the result $\brho_\star$ of the $\ell_1$ minimization, which
is no longer guaranteed to be unique.  We show that there are cases
where the minimization may be useful.  Specifically, we prove that
when the unknown sources/scatterers are located at points or clusters
of points that are sufficiently well separated, an $\ell_1$ minimizer
$\brho_\star$ is supported in the vicinity of these points. While the
entries of $\brho_\star$ may not be close in the point-wise sense to
those of $\brho$, their average over such vicinities are close to the
true values in $\brho$ in the case of well separated points, or the
averages of the true values in the case of clusters of points. That is
to say, $\ell_1$ optimization gives an effective vector of
source/scatterer amplitudes averaged locally around the points in its
support.  Note that question (2) was also investigated in
\cite{fannjiang2012coherence}, where novel algorithms for imaging well
separated sources have been introduced and analyzed. Our study
complements the results in \cite{fannjiang2012coherence} by analyzing
directly the performance of the $\ell_1$ minimization and $\ell_1$-penalty,  and also
considering clusters of sources/scatterers.

The paper is organized as follows. In section \ref{sect:formulate} we
formulate the problem, introduce notation, and describe the relation
between imaging sources vs. scatterers. Question (1) is studied in
section \ref{sect:parax}.  We describe the paraxial scaling regime and
derive resolution bounds that depend on the sparsity $s$ of the
imaging scene.  In section \ref{sect:cont} we study question (2).  In
both sections we begin with the statement of results and numerical
illustrations, and end with the proofs.  A summary is in section
\ref{sect:sum}.

\section{Formulation of the imaging problem}
\label{sect:formulate}
We formulate first the basic problem of imaging $s$ point-like sources
with a remote array of sensors that record the incoming sound waves.
The generalization to the inverse scattering problem is described in
section \ref{sect:F1.2}.

Suppose that there are $s$ unknown sources located at points $\vy_j$
in the imaging region $W \subset \real^3$, emiting signals $\hat
f_j(\om)$ at frequency $\om$, for $j = 1, \ldots, s$. The hat stands
for Fourier transform with respect to time, and reminds us that we
work in the frequency domain.  The receivers are at locations $\vx_r
\in \cA$, for $r = 1, \ldots, M_r$, where $\cA$ is a set on the
measuring surface, called the array aperture. The sound pressure wave
measured at $\vx_r$ and frequency $\om$ is modeled by
\begin{equation}
\hat p(\om,\vx_r) = \sum_{j=1}^s \hat f_j(\om) \hat
G(\om,\vx_r,\vy_j),
\label{eq:F1}
\end{equation}
where $\hat G$ is the outgoing Green's function of Helmholtz's
equation. The propagation is through a homogeneous medium with sound
speed $c$, and the Green's function is
\begin{equation}
\hat G(\om,\vx_r,\vy_j) = \frac{e^{i k |\vx_r-\vy_j|}}{4 \pi |\vx_r-\vy_j|},
\label{eq:F2}
\end{equation}
where $k = \om/c$ is the wavenumber.  The inverse source problem is to
determine $\{\hat f_j\}_{j = 1, \ldots, s}$ from the measurements
(\ref{eq:F1}) at one or more frequencies $\om$.

\subsection{Imaging sources with $\ell_1$ optimization}
\label{sect:F1.1}
To state the inverse problem as an $\ell_1$ optimization we discretize
$W$ with a regular grid of rectangular prisms, and let $W_N \subset
\real^3$ be the set of $N$ grid points denoted by $\vz_j$.  The
lengths of the edges of the rectangular prisms are the components of
the vector $\vbh \in \real^3$, called the mesh size. The sources may
be on or off the grid. If they are on the grid, as assumed in section
\ref{sect:parax}, we denote by $\cS$ the set of indexes of the grid
points that support them.  Explicitly, we define the bijective map
$J:\{1, \ldots, s\} \to \cS \subset \{1 \ldots, N\}$, such that
$\vz_{_{J(j)}} = \vy_j$, for $j = 1, \ldots, s$, and $\cS = \{J(1),
\ldots, J(s)\}$. When the sources are not on the grid, we let $\cS$
index the nearest points $\vz_j$ to each source, as explained in more
detail in section \ref{sect:cont}.  We assume henceforth that $N \gg
s$, meaning that the imaging scene is sparse.

Let $\bd\in \complex^M$ be the data vector, with components $\hat
p(\om_l,\vx_r)$, for $l = 1, \ldots, M_\om$ and $r = 1,\ldots,
M_r$. The number of measurements is $M = M_\om M_r$. Let also $\bG \in
\complex^{M\times N}$ be the sensing matrix, with entries defined by
$\hat G(\om_l,\vx_r,\vz_j)/\alpha_j$, where $\alpha_j$ normalizes the
columns of $\bG$, denoted by $\bg_j$.  Absorbing the normalization
constants in the vector $\brho$ of unknowns, we obtain the linear
system
\begin{equation}
\bG \brho = \bd.
\label{eq:F3}
\end{equation}
For single frequency measurements at $\om = \om_1$ the normalization
constant is
\[
\alpha_j = \Big(\sum_{r=1}^{M_r} |\hat
  G(\om,\vx_r,\vz_j)|^2\Big)^{1/2} = \Big(\sum_{r=1}^{M_r}
  \frac{1}{16 \pi^2 |\vx_r-\vz_j|^2}\Big)^{1/2},
\]
so that when the sources lie on the grid there are $s$ non-zero
entries in $\brho$, equal to
\begin{equation}
\rho_j = \alpha_{J^{-1}(j)} \hat f_{J^{-1}(j)}(\om), \qquad j \in \cS.
\label{eq:F4}
\end{equation}
Here $J^{-1}:\cS \to \{1, \ldots, s\}$ is the inverse of the mapping
$J$.  For multiple frequency measurements we simplify the problem by
letting
\begin{equation}
\hat f_j(\om) = \hat f(\om) R_j, \qquad \forall j = 1, \ldots, s,
\label{eq:F5}
\end{equation} 
so that all the sources emit the same known signal $\hat f(\om)$
multiplied by an unknown complex amplitude $R_j$. This simplification
is motivated by the inverse scattering problem described in section
\ref{sect:F1.2}. It keeps the same number $N$ of unknowns as in the
single frequency case, although we have more measurements.  The
normalization constants are
\[
\alpha_j = \|\hat f\|_{2} \Big(\sum_{r=1}^{M_r} \frac{1}{16 \pi^2
    |\vx_r-\vy_j|^2}\Big)^{1/2}, \qquad \|\hat f\|_2 =
\Big(\sum_{l=1}^{M_\om} |\hat f(\om_l)|^2\Big)^{1/2},
\] 
and the non-zero entries of $\brho$ equal
\begin{equation}
 \rho_j = \alpha_{J^{-1}(j)} R_{J^{-1}(j)}, \qquad j \in \cS.
\label{eq:F6}
\end{equation}
Note that we could have written the multiple frequency problem for
$M_\om N$ unknowns, the Fourier coefficients $\hat f_j(\om_l)$ of the
signals emitted by the sources. However, at each frequency these have
the same spatial support, so another optimization approach, known as
Multiple Measurement Vector (MMV) \cite{cotter2005sparse} would be
more appropriate. For the purpose of this paper it suffices to
consider the simpler model (\ref{eq:F6}).

The $\ell_1$ optimization (basis pursuit) formulation of the inverse
source problem is
\begin{equation}
\min_{\brho \in \complex^N} \|\brho\|_{1} \quad \mbox{such that} \quad
\bG \brho = \bd,
\label{eq:F7}
\end{equation}
where $\|\brho\|_1 = \sum_{j=1}^N |\rho_j|.$ Our goal in section
\ref{sect:parax} is to determine bounds on the
mesh size $\vbh$ so that (\ref{eq:F7}) has a unique $s$ sparse
solution, equal to the true $\brho$ defined in (\ref{eq:F4}) and
(\ref{eq:F6}). The analysis is based on the next lemma, following from
\cite{tropp2004greed,tropp2006just,DonElad}.

\vspace{0.1in}
\begin{lemma}
\label{lem.1}
Suppose that (\ref{eq:F3}) has an $s$ sparse solution $\brho$ and that
the cumulative mutual coherence $\mu(\bG,s)$ of matrix $\bG$ with
columns $\bg_j$ of Euclidian length equal to one satisfies
\begin{equation}
\mu(\bG,s) =  \max_{j=1, \ldots, N}\max_{|S| = s} \sum_{q\in S, q
  \ne j} |\lin \bg_q, \bg_j\rin| < \frac{1}{2},
\label{eq:F8}
\end{equation}
where $\lin \cdot,\cdot \rin$ denotes the usual inner product in
$\complex^N$, and $S$ is a set of cardinality $|S|$.  Then $\brho$ is
the unique $s$ sparse solution of (\ref{eq:F3}) and the unique
minimizer of (\ref{eq:F7}).
\end{lemma}

To deal with noise and modeling (discretization) error, we also
consider in section \ref{sect:cont} the $\ell_1-$penalty problem
\cite{tropp2006just}
\begin{equation}
  \min_{\brho \in \complex^N} \sL(\brho), \qquad \sL(\brho) =
  \frac{1}{2} \|\bG \brho - \bd\|_2^2 + \gamma \|\brho\|_1,
  \label{eq:F7b}
\end{equation}
with parameter $\gamma$ accounting for the trade-off between the
approximation error and the sparsity of the unknown vector $\brho$.

\subsection{Imaging point scatterers}
\label{sect:F1.2}
The problem of imaging a sparse scene of point scatterers at $\vy_j$
for $j = 1, \ldots, s,$ can be written as one of imaging $s$ sources,
as we now explain.

Let the array probe the medium with a signal $\hat f(\om)$ emitted
from the sensor at $\vx_e$. Using the Foldy-Lax model
\cite{foldy1945multiple,lax1951multiple} we write the scattered wave
at $\vx_r$ in the form (\ref{eq:F1}), with effective sources at
$\{\vy_j\}_{j = 1, \ldots, s}$ emitting signals
\begin{equation}
\hat f_j(\om) = \hat f(\om) R_j \hat u_j(\om).
\label{eq:F9}
\end{equation}
Here $R_j$ is the reflectivity of the $j-$th scatterer, and $\hat u_j$
is the wave that illuminates it. It is given by the sum of the
incident wave $\hat G(\om,\vx_e,\vy_j)$ and the wave scattered at the
other points
\begin{equation}
\hat u_j(\om) = \hat G(\om,\vx_e,\vy_j) + \sum_{l = 1}^s
(1-\delta_{lj})R_l \hat G(\om,\vy_l,\vy_j) \hat u_l(\om), \qquad
\forall j = 1, \ldots, s,
\label{eq:F10}
\end{equation}
where $\delta_{lj}$ is the Kronecker delta.

In the Born approximation we neglect the sum in (\ref{eq:F10}), and
simplify (\ref{eq:F9}) as
\begin{equation}
\hat f_j(\om) = \hat f(\om) R_j \hat G(\om,\vx_e,\vy_j).
\label{eq:F11}
\end{equation}
This can be written in the form (\ref{eq:F3}) with entries of $\brho$
like in (\ref{eq:F6}), and slightly redefined matrix $\bG$ and
normalization constant
\[
\alpha_j = \|\hat f\|_{2} \sqrt{\sum_{r=1}^{M_r} \frac{1}{(4 \pi)^4
    |\vx_r-\vy_j|^2 |\vx_e-\vy_j|^2}}.
\]

Multiple scattering effects can be included by solving the Foldy-Lax
equations (\ref{eq:F10}) or, equivalently, the linear system
\[
{\bf Q}\hat {\bf u} = \left( \begin{matrix} \hat G(\om,\vx_e,\vy_1)
  \\ \vdots \\ \hat G(\om,\vx_e,\vy_s) \end{matrix} \right),
\]
with $\hat {\bf u}$ the vector with components $\hat u_j$, and ${\bf
  Q} = (Q_{jl})_{j,l = 1, \ldots s}$ the matrix with entries \[ Q_{jl}
= \delta_{lj} - (1-\delta_{jl}) \hat G(\om,\vy_l,\vy_j) R_l.\] Note
that ${\bf Q}$ is a perturbation of the $s \times s$ identity matrix,
and depending on the magnitude of the reflectivities and the distance
between the scatterers, it is invertible. Again we can write the
problem in the form (\ref{eq:F3}) with entries of $\brho$ like in
(\ref{eq:F6}), except that now $\alpha_j$ 
 are more complicated and
depend on the unknown reflectivity. An elegant solution of this
nonlinear problem is in \cite{chai2014imaging}.  It amounts to solving
a source imaging problem like (\ref{eq:F7}), to determine the
locations $\vy_j$, for $j = 1, \ldots, s.$ Because there are multiple
emitters the authors use an MMV approach.  Then the reflectivities are
estimated using the Foldy-Lax model.

Given the relation between inverse scattering and source problems
describe above, we focus attention henceforth on imaging sparse scenes
of sources using (\ref{eq:F7}) or (\ref{eq:F7b}).
\section{Imaging with small arrays}
\label{sect:parax}
\begin{figure}[t]
\centerline{\includegraphics[width = 0.7
    \textwidth]{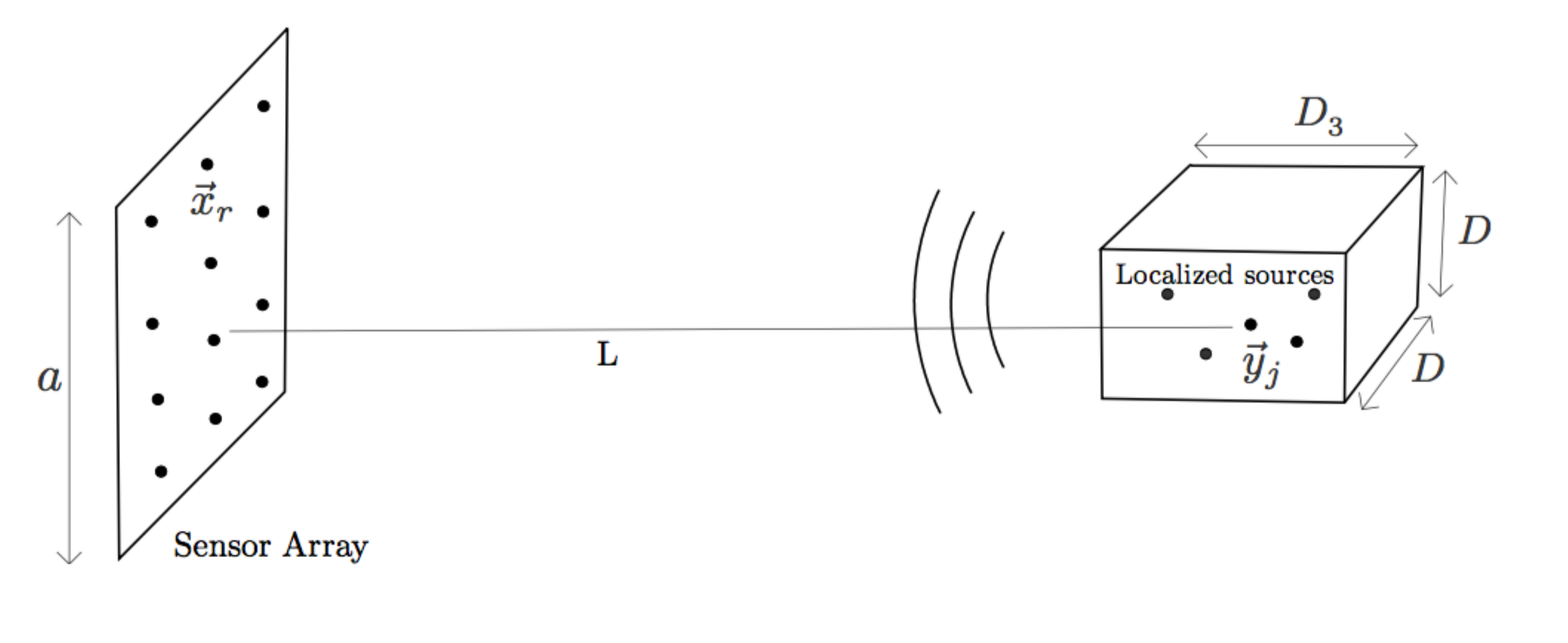}}

\vspace{-0.2in}
\caption{Schematic for the paraxial setup. }
\label{fig:parax}
\end{figure}

The setup is illustrated in Figure \ref{fig:parax}. We consider a
planar square array of aperture size $a$, and a coordinate system with
origin at the center of the array and range axis orthogonal to it. The
locations of the receivers are $\vx_r = (\bx_r,0)$, with $\bx_r =
(x_{1,r},x_{2,r})$ and $|x_{1,r}|, |x_{2,r}| \le a/2$, for $r = 1,
\ldots M_r$. The imaging region $W$ is a rectangular prism with center
on the range axis, at distance $L$ from the array. It has a square
side $[-D/2,D/2] \times [-D/2,D/2]$ in the cross-range plane, parallel
to the array, and length $D_3$ in the range direction. The
discretization of $W$ has grid points $\vz_j = (\bz_j,z_{3,j})$, with
cross-range vector $\bz_j = (z_{1,j},z_{2,j})$ and range $z_{3,j}$,
and the mesh size is $\vbh = (h,h,h_3)$. Our goal in this section
is to estimate $\vbh$ so that the $\ell_1$ optimization problem
(\ref{eq:F7}) determines exactly the unknown sources supported at
$\vz_j$ for $j \in \cS$, a set of cardinality $s$.

We begin in section \ref{sect:parax_scaling} with the scaling regime
and the paraxial model of wave propagation. The resolution limits are
stated and illustrated with numerical simulations in sections
\ref{sect:parax_sf} and \ref{sect:parax_bb}.  The setup of the
simulations is described in appendix \ref{ap:numeric}. The proofs are
in section \ref{sect:parax_proofs}.

\subsection{Scaling regime and the paraxial model}
\label{sect:parax_scaling}
The scaling regime is defined by the relation between the important
length scales: the wavelength $\la$, the range scale $L$, the aperture
$a$, and the size $D$ and $D_3$ of the imaging region.  We assume for
now a single frequency $\om$, and refer to section \ref{sect:parax_bb}
for the multi frequency case where another important scale arises, the
bandwidth $B$.

The scales are ordered as 
\begin{equation}
\la \ll D \ll a \ll L, \qquad D_3 \ll L,
\label{eq:P1}
\end{equation}
and satisfy the following assumptions
\begin{align}
\frac{a^2}{\la L} &\gtrsim 1, \label{eq:P2} \qquad \frac{a^2}{\la L}
\frac{D_3}{L}\gtrsim 1, \\ \frac{D^2}{\la L} & \ll 1, \label{eq:P4}
\quad \frac{a^2}{\la L} \frac{a D}{L^2} \ll 1, \quad \frac{a^2}{\la L}
\left(\frac{D_3}{L}\right)^2  \ll 1, \quad  \frac{a^2}{\la L}
\left(\frac{a}{L}\right)^2 \frac{D_3}{L}  \ll 1. 
\end{align}
Roughly, conditions (\ref{eq:P4}) say that the imaging region is small
enough and far enough from the array, so that we can linearize phases
of the Green's functions in $\vz_j$, for all $j = 1, \ldots, N.$
Physically, this means that when viewed from the imaging region, the
wave fronts appear planar. The array aperture $a$ is small with
respect to the distance $L$ of propagation of the waves, but equations
(\ref{eq:P2}) say that the Fresnel number is large, so we have
diffraction effects.

We show in appendix \ref{ap:parax_deriv} that 
\begin{equation}
\label{eq:PGG}
\sum_{r=1}^{M_r} \hat G(\om,\vx_r,\vz_j) \overline{\hat
  G(\om,\vx_r,\vz_q)} \approx \frac{e^{i k (z_{3,j}-z_{3,q})}}{(4 \pi
  L)^2} \sum_{r=1}^{M_r}e^{-i k\big[ \frac{|\bx_r|^2(z_{3,j}-z_{3,q})}{2 L^2}
     + \frac{\bx_r \cdot (\bz_j-\bz_q)}{L}\big]},
\end{equation}
where the bar denotes complex conjugate. 

\begin{remark} \label{rem:z3}
The terms $e^{i k z_{3}}$   in  (\ref{eq:PGG}) are highly oscillatory when $k$ is large,
but can be absorbed in the vector of unknowns. This is convenient because it implies that when two points $\vz_j$ and $\vz_q$ are close to each other,
the inner product of the $j-$th and $q-$th columns of the scaled sensing matrix is close to one. 
Let $\mathscr{Z}$ be the diagonal matrix  with $j^\text{th}$ entry $e^{i k z_{3,j}}$, for $j = 1, \ldots, N$,
and  rewrite the linear system  (\ref{eq:F3}) as $\bG \mathscr{Z} (\mathscr{Z}^{-1} \brho)  = \bd$.
To
 simplify the presentation we denote henceforth by $\bG$ the new sensing matrix 
 $\bG \mathscr{Z}$  and by
$\brho$ the new unknown vector $\mathscr{Z}^{-1} \brho$, so that the scaled system looks  the 
same as (\ref{eq:F3}). 
The Green's function with the large phase removed is still denoted by 
$\hat G$ , and satisfies
\begin{equation}
\sum_{r=1}^{M_r} \hat G(\om,\vx_r,\vz_j) \overline{\hat
  G(\om,\vx_r,\vz_q)} \approx \frac{1}{(4 \pi
  L)^2} \sum_{r=1}^{M_r}e^{-i k\big[ \frac{|\bx_r|^2(z_{3,j}-z_{3,q})}{2 L^2}
     + \frac{\bx_r \cdot (\bz_j-\bz_q)}{L}\big]}.
\label{eq:P8}
\end{equation}

\end{remark}

 Assuming that the receivers
are on a square grid of spacing $h_{_\cA}$,
satisfying the scaling relations
\begin{equation}
h_{_\cA} \ll \frac{ \la L}{D}, \qquad h_{_\cA} \ll \frac{\la L^2}{a D_3},
\label{eq:P9}
\end{equation}
we see that the exponential in (\ref{eq:P8}) is approximately constant
in each grid cell in $\cA$. This allows us to use the continuous
aperture approximation in the analysis, where the sum over $r$ can be
replaced by the integral over $\cA = [-a/2,a/2] \times [-a/2,a/2]$,
\begin{equation}
\sum_{r=1}^{M_r} \hat G(\om,\vx_r,\vz_j) \overline{\hat
  G(\om,\vx_r,\vz_q)} \approx \frac{1}{(4 \pi
  L h_{_\cA})^2} \hspace{-0.05in}\int_{\cA} d \bx \, e^{- i
  k\big[\frac{ |\bx|^2(z_{3,j}-z_{3,q})}{2 L^2} + \frac{\bx \cdot
      (\bz_j-\bz_q)}{L}\big]}.
\label{eq:P10}
\end{equation}
With our discretization the number of unknowns is $ N =
\left({D}/{h}\right)^2 {D_3}/{h_3}$ and the problem is underdetermined
when the number of measurements $ M = M_r =
\left({a}/{h_{_\cA}}\right)^2 $ satisfies $M < N$ or, equivalently,
\begin{equation}
h_{_\cA} > a \frac{h}{D} \sqrt{\frac{h_3}{D_3}},
\label{eq:P13}
\end{equation}
which is consistent with (\ref{eq:P9}) for $h \ll D$ and $h_3 \ll D_3$.

\subsection{Single frequency resolution limits}
\label{sect:parax_sf}
Using the paraxial model for our sensing matrix $\bG$, we obtain from
(\ref{eq:P10}) and the definition of the cumulative coherence
$\mu(\bG,s)$ that
\begin{equation}
  \mu(\bG,s) = \hspace{-0.05in}\max_{j=1, \ldots, N} \max_{|\rm{S}| =
    s} \hspace{-0.05in}\sum_{q \in \rm{S}, q \ne j}\hspace{-0.05in}
  \cU\Big(\frac{z_{1,j}-z_{1,q}}{H},\frac{z_{3,j}-z_{3,q}}{H_3}\Big)
  \cU\Big(\frac{z_{2,j}-z_{2,q}}{H},\frac{z_{3,j}-z_{3,q}}{H_3}\Big),
  \label{eq:cumC}
\end{equation}
where 
\begin{equation}
  H = \frac{L}{ka} = \frac{\la L}{2 \pi a}, \qquad H_3 = \frac{2
    L^2}{k a^2} = \frac{\la L^2}{\pi a^2},
  \label{eq:P14}
\end{equation}
and $\cU(\beta,\eta)$ is the absolute value of the Fresnel integral
\begin{equation}
  \cU(\beta,\eta) = \left| \int_{-1/2}^{1/2} dt \, e^{-i \beta t - i
    \eta t^2}\right|.
  \label{eq:P15}
\end{equation}
The search set of cardinality $s$ is denoted by $\rm{S}$, to
distinguish it from the set of indexes of the true support points of
$\brho$, called calygraphic $\cS$.

As stated in Lemma \ref{lem.1}, unique recovery of a sparse $\brho$
with the $\ell_1$ minimization (\ref{eq:F7}) is guaranteed when
$\mu(\bG,s) < 1/2$. This criterion allows us to estimate the
resolution limit stated in the next two theorems.

\vspace{0.1in}
\begin{theorem}
  \label{thm.1}
  If the mesh size satisfies
  \begin{equation}
  h > h^\star = \frac{2}{\pi} \frac{\la L}{a}, \qquad h_3 > h^\star_3
  = \frac{16}{\pi} \frac{\la L^2}{a^2},
  \label{eq:P16}
  \end{equation}
  $\ell_1$ optimization recovers exactly two sources located at any
  distinct grid points.
\end{theorem}

\vspace{0.1in} \noindent We call the estimates $h^\star$ and
$h_3^\star$ the ``base resolution''.  They are the same, up to order
one constants, as the well known resolution limits in array imaging
\cite{born1999principles}.  We verified numerically the estimates
(\ref{eq:P16}) as follows. To check the value of $h^\star$, we solved
the optimization problem (\ref{eq:F7}), as explained in appendix
\ref{ap:numeric}, for data corresponding to a vector $\brho$ supported
at any two grid points offset in cross-range.  We determined the
smallest $h$ so that the relative error between $\brho$ and its
numerical reconstruction was less than $1\%$. A similar estimation was
done for $h^\star_3$, with $\brho$ supported at points offset in
range.  The results were close to those in Theorem \ref{thm.1}: $0.46
\la L/a$ for $h^\star$ and $3 \la L^2/a^2$ for $h^\star_3$.

\begin{figure}[t]
\begin{center}
\vspace{-0.8in}
\raisebox{-0.in}{{\includegraphics[width=0.42\textwidth]{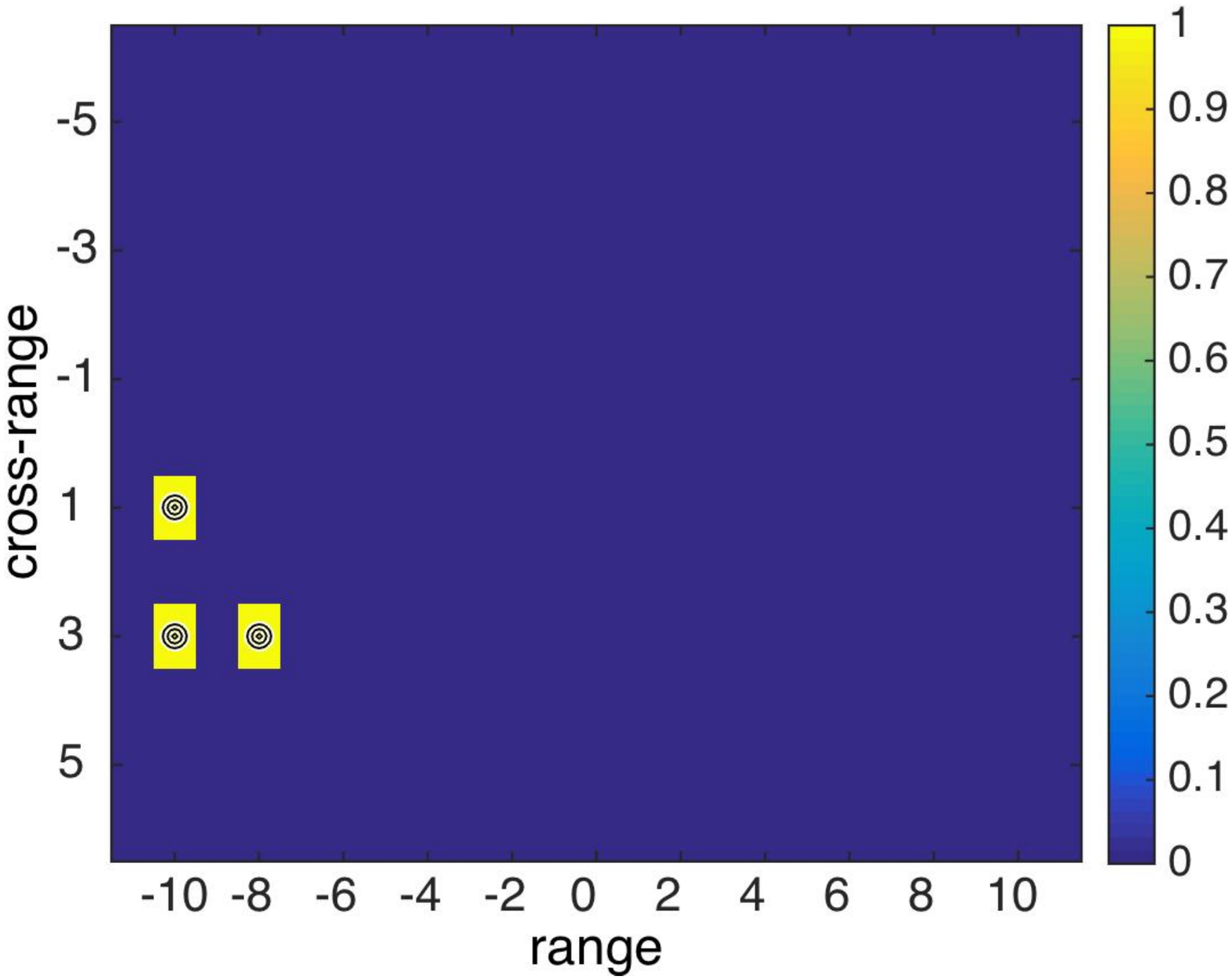}}}
\raisebox{-.84in}{
  {\includegraphics[width=0.48\textwidth]{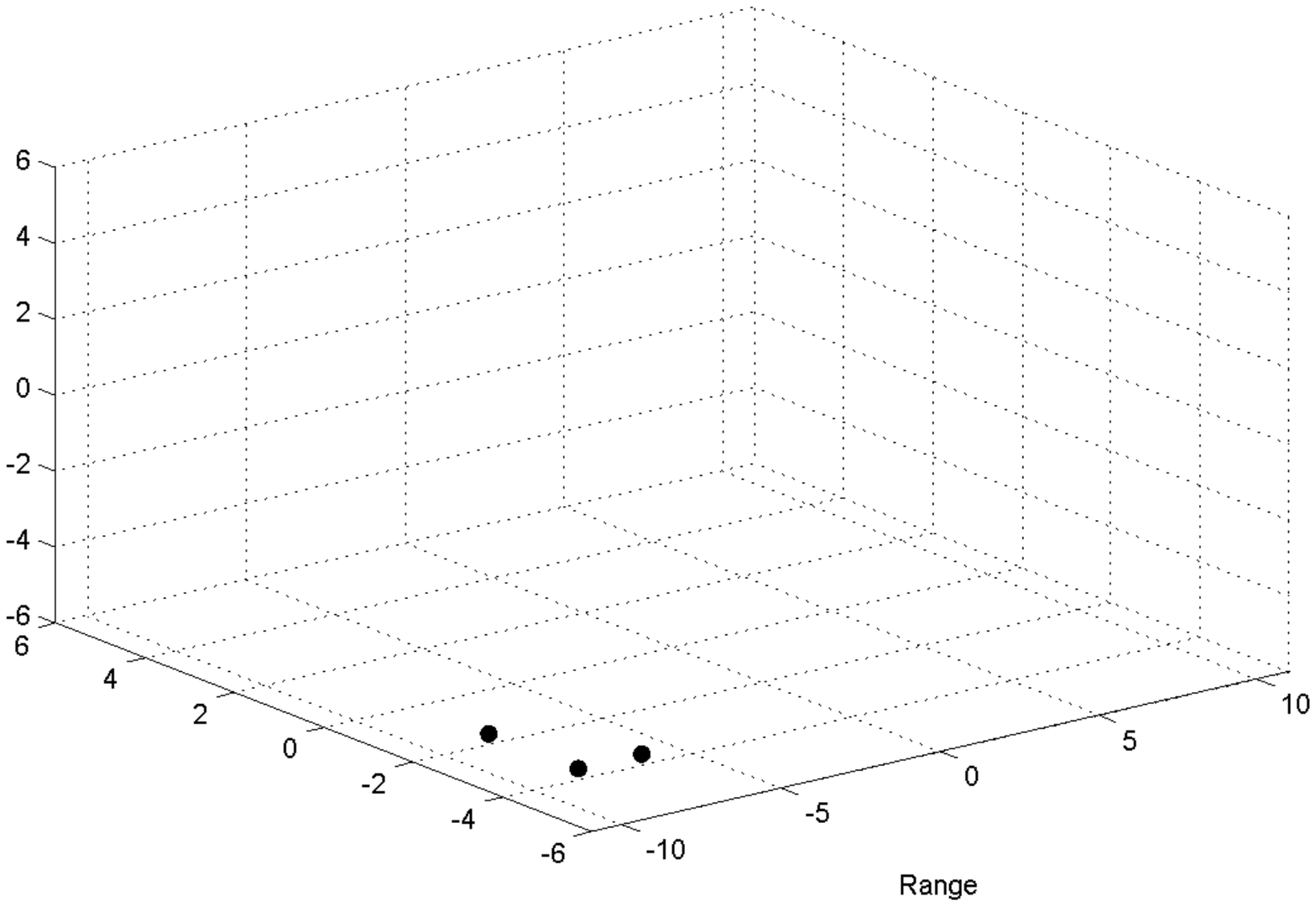}}}
                                                                                                  
\vspace{-0.88in} \raisebox{0
  in}{{\includegraphics[width=0.42\textwidth]{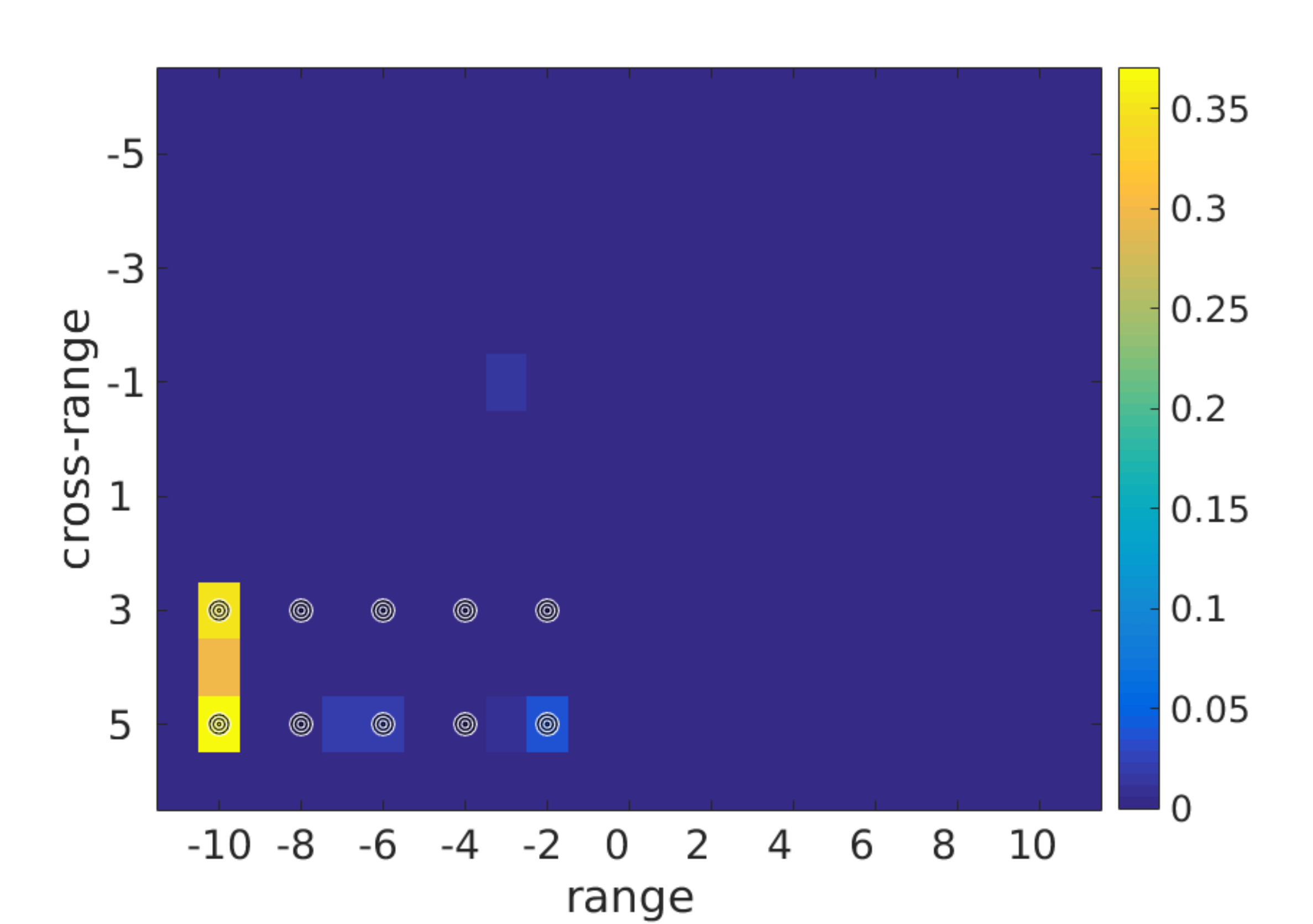}}}
\raisebox{0in}{
  {\includegraphics[width=0.48\textwidth]{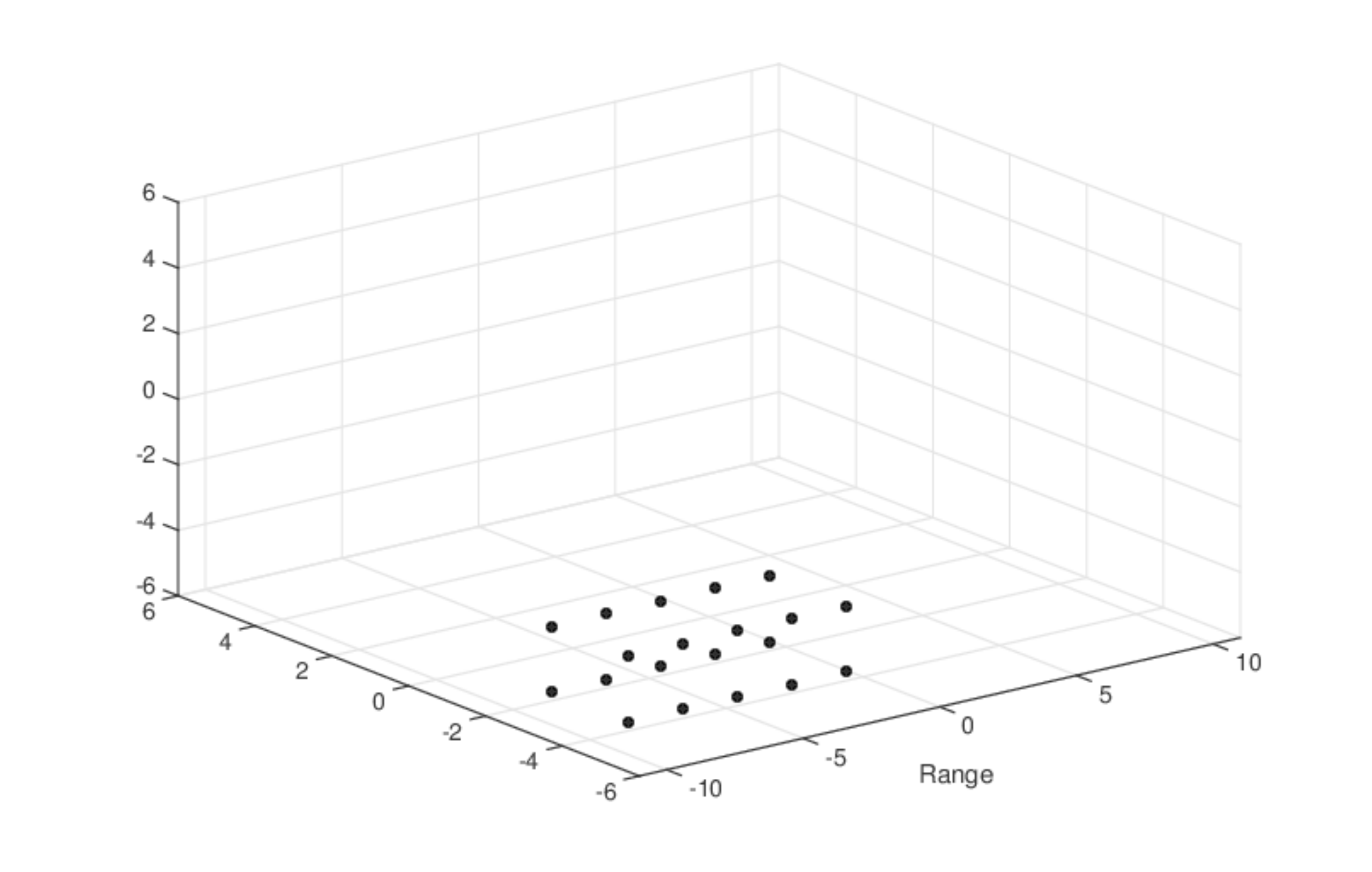}}}
                                                  

\end{center}

\vspace{-0.1in}\caption{Illustration of recovery of two imaging scenes
  discretized at base resolution $(h^\star,h^\star,h_{3}^\star)$. The
  results on the top line are for $3$ sources and on the bottom line
  for $20$ sources. The distribution of the sources is displayed in
  the right column. The left column shows a cross-section of the
  images. The range axis is in units of $h_3^\star$ and the
  cross-range axis in units of $h^\star$. The exact location of the
  sources is superposed on the images. They all have strength $\rho_j
  = 1$, for $j \in \cS$, and the magnitude of the reconstruction is
  shown with the color bar.  }
\label{fig:IllRes}
\end{figure}

The next result states that when there are more sources to estimate,
the resolution limits deteriorate. 
This is also illustrated in Figure
\ref{fig:IllRes}, where we display images discretized at base
resolution. The reconstruction is perfect for $3$ sources (top line),
but not for $20$ sources (bottom line).  How the resolution limits
deteriorate with $s$ depends on the distribution of the sources in the
imaging region. Our estimate in the next theorem accounts for worse
case scenarios.

\vspace{0.1in}
\begin{theorem}
  \label{thm.2}
  There exists a constant $C$ of order one such that if the mesh size
  satisfies the conditions 
  \begin{equation}
    \frac{h/h^\star}{h_3/h_3^\star} = 0(1) \quad \mbox{and}
    \quad \left[ \Big(\frac{h}{h^\star}\Big)^2 \frac{h_3}{h_3^\star}
      \right]^{1/3}> C s^{2/3},
    \label{eq:P17}
  \end{equation}
  $\ell_1$ optimization recovers exactly $s$ sources located at any
  distinct grid points.
\end{theorem}

\vspace{0.1in} \noindent The isotropic dilation of the mesh in
(\ref{eq:P17}) is for convenience, but the result generalizes to
anisotropic dilations, where the mesh is stretched much more in one
direction than the others. We refer to the proof in section
\ref{sect:parax_pf_sf} for details on the generalization. The
resolution decrease with $s$ predicted by Theorem \ref{thm.2} may be
traced to the slow decay with range offsets of the terms $|\lin
\bg_j,\bg_k \rin|$ summed in $\mu(\bG,s)$. This is also why
\[
h_3^\star/h^\star = 8 L/a \gg 1.
\]
Sources at different ranges may have strong interaction, hence they
must be further apart in order to get $\mu(\bG,s) < 1/2$. In the next
section we show that if the base range resolution improves, as it does
with broad band data, then there is almost no resolution loss with the
sparsity $s$.

\subsection{Broad band resolution limits}
\label{sect:parax_bb}
When we have $M_\om$ frequency measurements, the vector of unknowns is
defined as in (\ref{eq:F6}), and the rows of the sensing matrix $\bG$
are indexed by the receiver-frequency pair $(r,j) \in \{1, \ldots,
M_r\} \times \{1,\ldots, M_\om\}$. Let $\om_o$ be the central
frequency, so that
\begin{equation}
  |\om_j-\om_o| \lesssim B, \qquad j = 1, \ldots, M_\om,
  \label{eq:BB1}
\end{equation}
where $B$ is the bandwidth assumed to satisfy the scaling relation
\begin{equation}
 1 \lesssim \max \Big\{\frac{D}{\la_o L/a}, \frac{D_3}{\la_o L^2/a^2}
 \Big\} \ll \frac{\om_o}{B} \ll \Big(\frac{L}{a}\Big)^2.
  \label{eq:BB2}
\end{equation}
The lower bound says that $\om_o \gg B$, so that $\om_o$ is the scale
of all the measured frequencies $\om_j$, for $j = 1, \ldots
M_\om$. The upper bound implies $c/B \ll \la_o L^2/a^2$, where we
recall from the previous section that $\la_o L^2/a^2$ is, up to a
factor of order one, the base range resolution for single frequency
measurements at $\om = \om_o$. The next theorem states that the base
range resolution for multi frequency measurements is of order
$c/B$, so (\ref{eq:BB2}) implies a gain in range resolution.

Let us assume, for simplicity of the calculations, a Gaussian signal
\begin{equation}
  \hat f(\om_j) = e^{-\frac{(\om_j-\om_o)^2}{4 B^2}}.
  \label{eq:BB4}
\end{equation}
The results should extend to any signal with bandwidth $B$, with
modifications of the constants in the bounds. We show in appendix
\ref{ap:BBparax_deriv} that
\begin{align}
  \sum_{j=1}^{M_\om} \sum_{r=1}^{M_r} |\hat f(\om_j)|^2 \hat
  G(\om_j,\vx_r,\vz_q) \overline{\hat G(\om_j,\vx_r,\vz_{l})} \approx
  \frac{1}{(4 \pi L)^2} \sum_{j=1}^{M_\om}
  e^{-\frac{(\om_j-\om_o)^2}{2 B^2}+i \frac{(\om_j-\om_o)}{c}
      (z_{3,q} - z_{3,l})}\times \nonumber
    \\ \sum_{r=1}^{M_r} e^{- i k_o \Big[
        \frac{|\bx_r|^2 (z_{3,q}-z_{3,l})}{2L^2} + \frac{\bx_r \cdot
          (\bz_q-\bz_{l})}{L} \Big]},
  \label{eq:BB3}
\end{align}
where $k_o =\om_o/c$ is the central wavenumber, and we proceeded as in Remark \ref{rem:z3}
to absorb the large phases $e^{i k_o z_{3,q}}$ in the vector of unknowns. Assuming that the
array is discretized on a mesh with spacing $h_{_\cA}$ satisfying
(\ref{eq:P9}), we approximate the sum over $r$ by an integral over the
aperture $\cA$, as in the previous section. We also suppose that the
frequencies $\om_j$ are spaced at intervals $h_\om$ satisfying $ h_\om
\ll {c}/{D_3},$ so that we can write the sum over the frequencies as
an integral over the bandwidth. Equation (\ref{eq:BB3}) becomes
\begin{align}
\sum_{j=1}^{M_\om} \sum_{r=1}^{M_r} |\hat f(\om_j)|^2 \hat
G(\om_j,\vx_r,\vz_q) \overline{\hat G(\om_j,\vx_r,\vz_{l})} \approx
\frac{\sqrt{2 \pi} B}{(4 \pi L h_{_\cA})^2 h_\om} e^{-\frac{B^2(z_{3,q}-z_{3,l})^2}{2
      c^2}}\times \\ \int_{\cA} d \bx \, e^{- i k_o \Big[ \frac{|\bx|^2
          (z_{3,q}-z_{3,l})}{2L^2} + \frac{\bx \cdot
          (\bz_q-\bz_{l})}{L} \Big]}, \label{eq:BB6}
\end{align}  
and we can simplify it further by neglecting the quadratic phase in
the integral over the aperture. This is because 
\[
\frac{k_o |\bx|^2 |z_{3,q}-z_{3,l}|}{L^2} \lesssim O \Big(\frac{k_o
  a^2 c}{L^2 B}\Big) = O\Big(\frac{\om_o a^2}{B L^2}\Big) \ll 1,
\]
for range offsets in the support of order $c/B$ of the Gaussian factor
in (\ref{eq:BB6}). Integrating over the aperture and normalizing, we
arrive at the following model of the products of the columns of the
sensing matrix $\bG$,
\begin{equation}
  |\lin \bg_q,\bg_{l}\rin |= e^{-\frac{(z_{3,q}-z_{3,l})^2}{2
      (c/B)^2}} \Big| \mbox{sinc} \Big(\frac{z_{1,q}-z_{1,l}}{2
    L/(k_o a)}\Big) \Big| \Big| \mbox{sinc}
  \Big(\frac{z_{2,q}-z_{2,l}}{2L/(k_o a)}\Big) \Big|.
  \label{eq:BB7}
\end{equation}
These are the terms in the cumulative coherence $\mu(\bG,s)$, and the
resolution limits are as stated  next.

\vspace{0.05in}
\begin{theorem}
  \label{thm.3}
  Assume a sensing matrix $\bG$ with inner products of the columns
  defined by (\ref{eq:BB7}).  If the mesh size $\vbh = (h,h,h_3)$
  satisfies
  \begin{equation}
    h >  h^\star = \frac{2}{\pi} \frac{\la_o L}{a}, \qquad
    h_3 > h_3^\star = \sqrt{2 \ln 2} \frac{c}{B},
    \label{eq:BB8}
  \end{equation}
  $\ell_1$ optimization recovers exactly any two sources on the
  grid. Moreover, there exists an order one constant $C > 1$,
  independent of $s$, such that if
  \begin{equation}
    \frac{h}{h^\star}, \frac{h_3}{h_3^\star} > C \ln s,
    \label{eq:BB9}
  \end{equation}
  $\ell_1$ optimization recovers exactly any $s$ sources on the grid.
\end{theorem}

\vspace{0.05in} We call $h^\star$ and $h_3^\star$ the base resolution,
as in the previous section. While the cross-range resolution $h^\star$
is the same as in the single frequency case, the base range resolution
$h_3^\star$ is significantly better.  Moreover, there is little loss
of resolution at large $s$. The mesh size grows at most
logarithmically with $s$, as opposed to $s^{2/3}$ in the single
frequency case. This agrees with the known fact in array imaging that
bandwidth improves the quality of images.

\subsection{Proofs}
\label{sect:parax_proofs}
The proofs of Theorems \ref{thm.1} and \ref{thm.2} which estimate the
resolution in the single frequency case are in section
\ref{sect:parax_pf_sf}. The proof of the broad band result in Theorem
\ref{thm.3} is in section \ref{sect:parax_pf_sf}.

\subsubsection{Single frequency}
\label{sect:parax_pf_sf}
We begin with some basic bounds on the Fresnel integral
(\ref{eq:P16}). The simplest estimate is for $\eta = 0$, in
which case
\begin{equation}
  \cU(\beta,0) = \Big| \mbox{sinc}\big(\beta/2\big)\Big| \le \min
  \{1,{2}/{\beta}\}.
  \label{eq:PP1}
\end{equation}
For $\eta \ne 0$ we can change variables and rewrite (\ref{eq:P16})
in the form
\begin{align}
\cU(\beta,\eta) &= \frac{1}{\sqrt{\eta}} \left|
\int_{\frac{\beta-\eta}{2 \sqrt{\eta}}}^{\frac{\beta+\eta}{2
    \sqrt{\eta}}} dt \, e^{-i t^2} \right| \le \frac{1}{\sqrt{\eta}}
\left[ \left| \int_0^{\frac{\beta+\eta}{2 \sqrt{\eta}}} dt \Big(\cos
  t^2 - i \sin t^2\Big) \right| + \left| \int_0^{\frac{\beta-\eta}{2
      \sqrt{\eta}}} dt \Big(\cos t^2 - i \sin t^2\Big)\right|
  \right]\nonumber \\&\le \frac{2 \sqrt{2}}{\sqrt{\eta}},
  \label{eq:PP2}
\end{align}
for any $\eta \ne 0$ and $\beta \ge 0$, where we used that
\[
\Big|\int_0^\alpha dt \, \cos t^2 \Big| \le 1, \qquad
\Big|\int_0^\alpha dt \, \sin t^2 \Big| \le 1, \qquad \forall \,
\alpha \in \real.
\]
The final estimate  
\begin{equation}
  \cU(\beta,\eta) = \frac{1}{\sqrt{\eta}} \left|
  \int_{\frac{\beta-\eta}{2 \sqrt{\eta}}}^{\frac{\beta+\eta}{2
      \sqrt{\eta}}} dt \, e^{ i t^2} \right| \le \frac{\pi + 1
  }{\alpha}, \qquad \mbox{for}~ ~ \beta > \alpha + \eta, \quad \forall
  \, \alpha, \eta > 0,
  \label{eq:PP10}
\end{equation}
follows from contour integration, as shown in appendix
\ref{ap:contour}.

\begin{figure}[t]
\begin{center}
\includegraphics[width =
  0.55\textwidth]{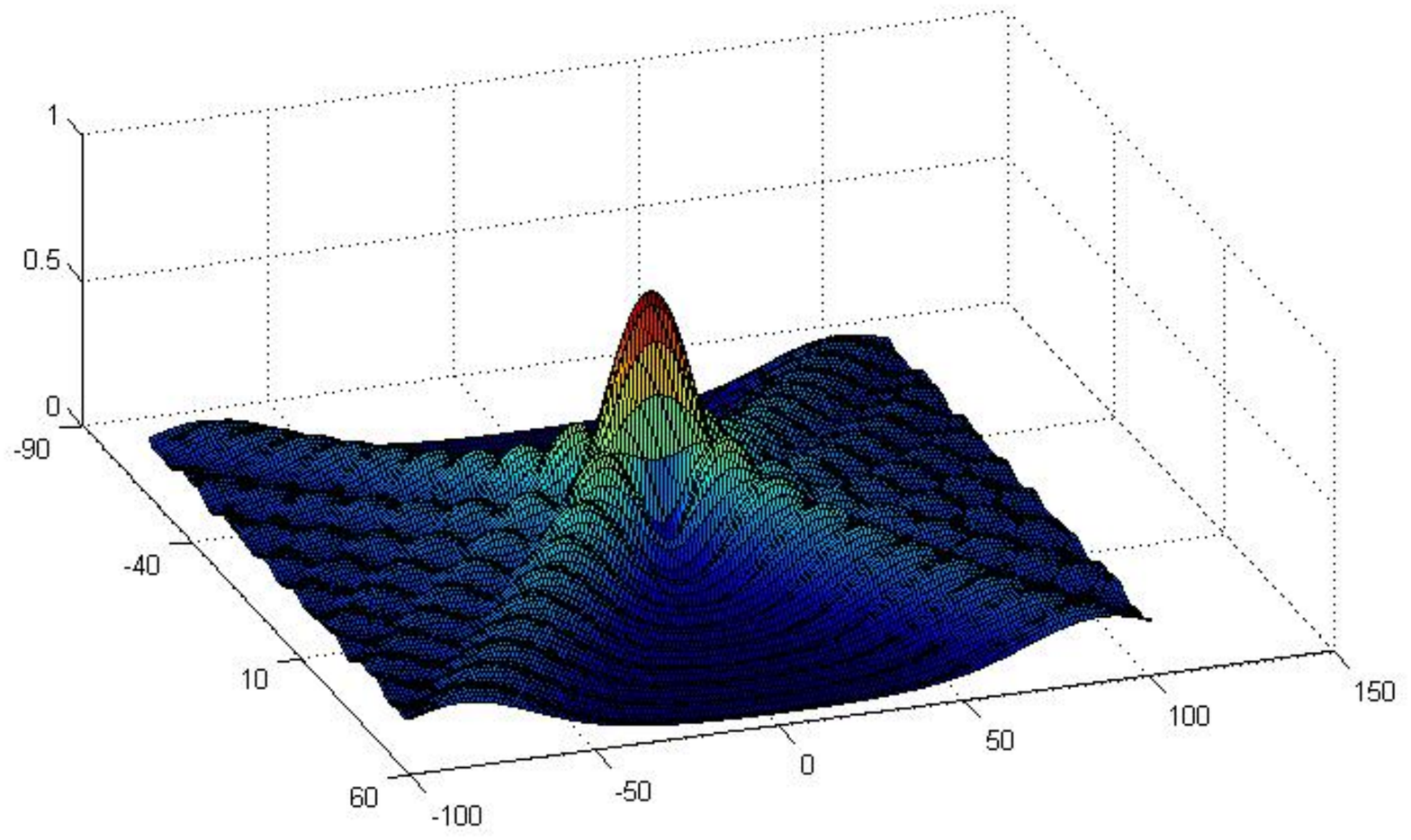}
\raisebox{0.3in}{\includegraphics[width =
    0.4\textwidth]{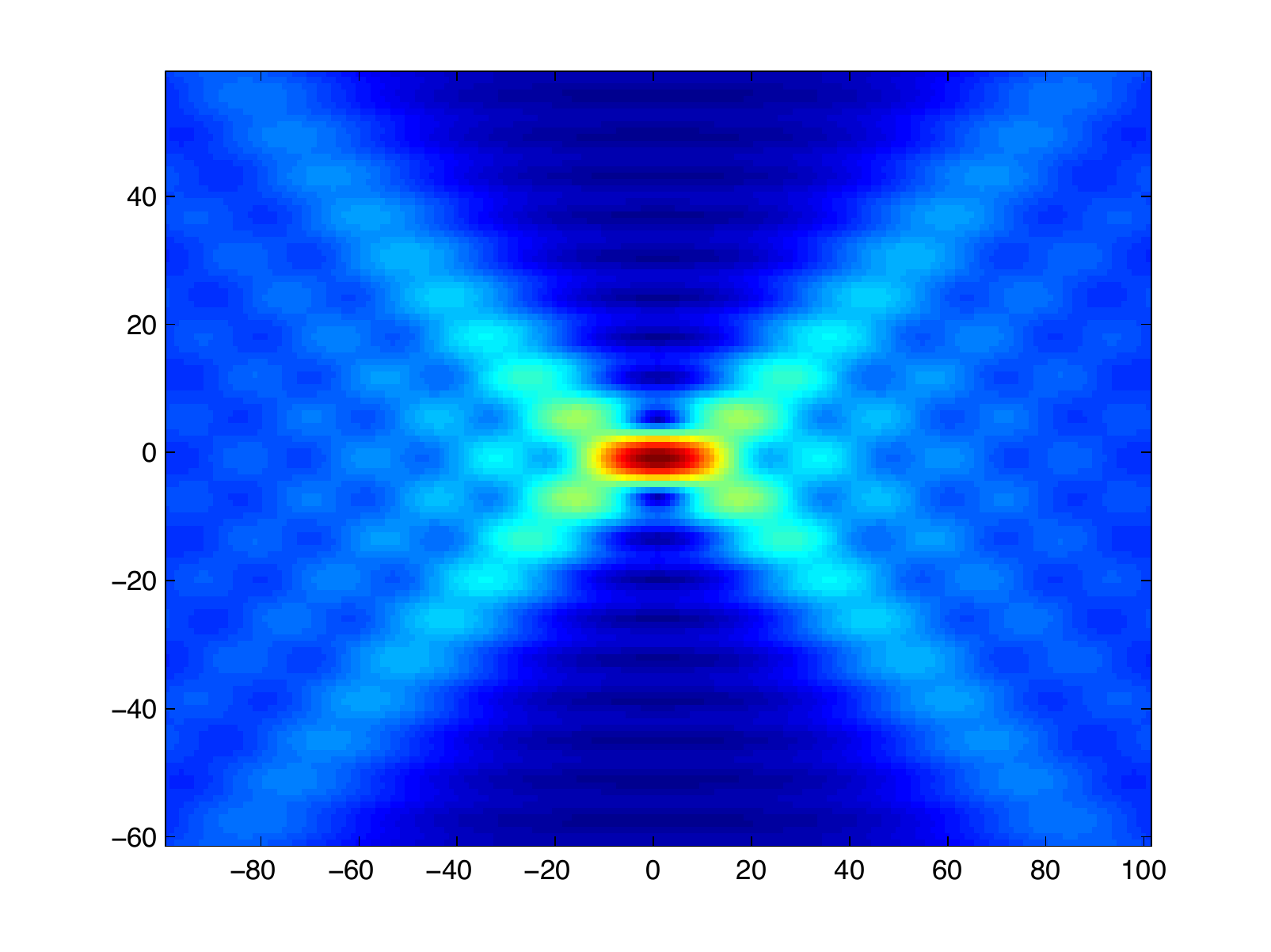}}
\end{center}

\vspace{-0.3in}\caption{Surface and top view display of the Fresnel
  integral $\cU(\beta,\eta)$ for $|\beta|,|\eta| \le 60$. In the right
  plot abscissa is $\eta$ and the ordinate is $\beta$.}
\label{fig:Fresnel}
\end{figure}

\textbf{Proof of Theorem \ref{thm.1}:} We wish to estimate $h$ and
$h_3$ so that $\mu(\bG,2) < 1/2$. The cumulative coherence for $s = 2$
is the same as the mutual coherence \cite{DonElad,tropp2006just}, and
in our case it takes the simple form
\begin{equation}
  \mu(\bG,2) = \max_{\vzeta \in \mathbb{Z}^3, \vzeta \ne 0}
  \cU\Big(\frac{h \zeta_1}{H},\frac{h_3 \zeta_3}{H_3}\Big)
  \cU\Big(\frac{h \zeta_2}{H},\frac{h_3 \zeta_3}{H_3}\Big).
  \label{eq:PP3}
\end{equation}
Here we used that the sources are on the grid and denote by $\vzeta =
(\zeta_1,\zeta_2,\zeta_3)$ vectors with integer components.  We
display in Figure \ref{fig:Fresnel} the Fresnel integral
$\cU(\beta,\eta)$, and note that it is bounded above by $1$, and
larger than $1/2$ for $\vzeta$ near the origin. When $\zeta_3 = 0$ we
get from (\ref{eq:PP1}) that 
\[
 \cU\Big(\frac{h \zeta_1}{H},0\Big) \cU\Big(\frac{h \zeta_2}{H},0
 \Big) \le \frac{2 H}{h},
 \]
because at least one of $\zeta_1$ and $\zeta_2$ is not equal to zero.
If $\zeta_3 \ne 0$ we have by (\ref{eq:PP2})
\[
\cU\Big(\frac{h \zeta_1}{H},\frac{h_3 \zeta_3}{H_3}\Big)
\cU\Big(\frac{h \zeta_2}{H},\frac{h_3 \zeta_3}{H_3}\Big) \le \frac{8
  H_3}{h_3}, \qquad \forall \, \zeta_1, \zeta_2 \in \mathbb{Z}.
\]
Thus, $\mu(\bG,2)$ is guaranteed to be less than $1/2$ if $h > 4 H =
h^\star$ and $h_3 > 16 H_3 = h^\star_3$. This concludes the proof of
Theorem \ref{thm.1}. $\Box$

\textbf{Proof of Theorem \ref{thm.2}:} Note from the expression
(\ref{eq:cumC}) of the cumulative coherence that it is translation
invariant in $\real^3$. Thus, we can fix the origin at one source
location and rewrite (\ref{eq:cumC}) as
\begin{equation}
  \mu(\bG,s) = \max_{|\Lambda| = s-1} \sum_{\vzeta \in \Lambda}
  \cU\Big(\frac{h \zeta_{1}}{H},\frac{h_3 \zeta_{3}}{H_3}\Big)
  \cU\Big(\frac{h \zeta_{2}}{H},\frac{h_3 \zeta_{3}}{H_3}\Big),
  \label{eq:PP4}
\end{equation}
where $\Lambda \subset \mathbb{Z}^3 \setminus\{0\}$ is a set of
cardinality $s-1$. The proof of the theorem follows from the bound on
$\mu(\bG,s)$ stated in the next lemma and the relation $h^\star = 4 H$
and $h_3^\star = 16 H_3$ established above. The constant $C$ in the
lemma is the same as in (\ref{eq:P17}).

\vspace{0.1in}
\begin{lemma}
  \label{lem.2}
  There exists constants $C$, $C_1$ and $C_2$ of order one such that
  \begin{align}
    \max_{|\Lambda| = s-1} \sum_{\vzeta \in \Lambda} \cU\Big(\frac{h
      \zeta_{1}}{H},\frac{h_3 \zeta_{3}}{H_3}\Big) \cU\Big(\frac{h
      \zeta_{2}}{H},\frac{h_3 \zeta_{3}}{H_3}\Big) \le 2^{5/3} C
    \left[\frac{s^2}{(h/H)^2 h_3/H_3}\right]^{1/3} + \nonumber \\ C_1
    \frac{H_3}{h_3} \ln s + C_2 \Big(\frac{H}{h}\Big)^2
    \ln^2 \hspace{-0.05in}s.
    \label{eq:lem2}
  \end{align}
\end{lemma}
The assumption in Theorem \ref{thm.2} that the mesh is dilated in an
isotropic fashion, and the relations $H = h^\star/4$ and $H_3 =
h_3^\star/16$ established above imply
\[
\frac{\left[\frac{s^2}{(h/H)^2 h_3/H_3}\right]^{1/3}}{\frac{H_3}{h_3}
  \ln s} = \frac{s^{2/3}}{\ln s} \left[\left(
  \frac{h_3}{H_3}\right)/\left(\frac{h}{H} \right)\right]^{2/3} =
  O\left(\frac{s^{2/3}}{\ln s}\right),
\]
and
\[
\frac{\left[\frac{s^2}{(h/H)^2
      h_3/H_3}\right]^{1/3}}{\Big(\frac{H}{h}\Big)^2
  \ln^2 \hspace{-0.05in} s} = \frac{s^{2/3}}{\ln^2 s} \frac{h}{H}
\left[\left( \frac{h}{H}\right)/ \left(\frac{h_3}{H_3}
  \right)\right]^{1/3} = O \left( \frac{s^{2/3}}{\ln^2\hspace{-0.05in}
  s} \frac{h}{H} \right).
\]
Thus, the first term in the right hand side of (\ref{eq:lem2})
dominates the others for large $s$ and $h \gtrsim h_\star = 4 H,$ and
the result stated in Theorem \ref{thm.2} follows. For anisotropic
dilations of the mesh the logarithmic terms in (\ref{eq:lem2}) may
become important. For example, when $h_3/H_3 > (h/H)^4
s^2/\ln^6\hspace{-0.05in}s$, the last term in (\ref{eq:lem2})
dominates the bound, and we can get a unique solution for a modest
mesh stretch in the cross-range direction $h/H = O(\ln s)$ at the
expense of a very large stretch in range $h_3/H_3 =
O(s^2/\ln^2\hspace{-0.05in}s)$.

\vspace{0.1in} \textbf{Proof of Lemma \ref{lem.2}:} Let
$\Lambda_\star$ be the set on which the maximum in (\ref{eq:PP4}) is
achieved. It is difficult to determine $\Lambda_\star$ explicitly, so
we construct another set $\Lambda_\tau$, which allows us to bound the
cumulative coherence. We use the behavior of the Fresnel integral
$\cU$, shown in Figure \ref{fig:Fresnel}, to guide us in the
construction. We write $\Lambda_\tau$ as the union of two sets
$\Lambda_\tau^+$ and $\Lambda_\tau^0$. The first set contains the left
and right cones (in Figure \ref{fig:Fresnel} they are defined by the
diagonal lines $|\beta| = |\eta|$), where $\cU$ displays a slower
decay, as well as a vicinity of the origin
\begin{equation}
   \Lambda_\tau^+ = \left\{ \vzeta \in \mathbb{Z}^3~ \mbox{s.t.}~ 0 <
   |\zeta_3| \le \frac{\tau}{h_3/H_3}, ~ ~ |\zeta_1|, |\zeta_2| \le
   \frac{H}{h} \left[\frac{\tau}{\sqrt{h_3 |\zeta_3|/H_3}} +
     \frac{h_3}{H_3} |\zeta_3| \right]\right\}.
  \label{eq:PP7}
\end{equation}
The second set is for the points with $\zeta_3 = 0$,
\begin{equation}
  \Lambda_\tau^0 = \left\{ \vzeta \in \mathbb{Z}^3\setminus \{0\} ~
  \mbox{s.t.}~ |\zeta_1|, |\zeta_2| \le \frac{\tau}{h/H}\right\}.
  \label{eq:PP6}
\end{equation}
The parameter $\tau > 1$ is used to control the volume of
$\Lambda_\tau = \Lambda_\tau^+ \cup \Lambda_\tau^0$, so that it
contains at least $s$ grid points,
\begin{equation}
\tau = \min \left\{ \left[\frac{3s}{8}
  \Big(\frac{h}{H}\Big)^2 \frac{h_3}{H_3}\right]^{1/3} ,
\frac{sh}{H}, \frac{sh_3}{H_3} \right\}.
\label{eq:PP5}
\end{equation}

The proof consists of two steps. First we derive the bound
\begin{align}
  \sum_{\zeta \in \Lambda_\star} \cU\Big(\frac{h
    \zeta_{1,q}}{H},\frac{h_3 \zeta_{3,q}}{H_3}\Big) \cU\Big(\frac{h
    \zeta_{2,q}}{H},\frac{h_3 \zeta_{3,q}}{H_3}\Big) \le \sum_{\zeta
    \in \Lambda_\tau} \mathcal{Z}(\vzeta),
  \label{eq:PP8}
\end{align}
where
\begin{equation}
  \mathcal{Z}(\vzeta) = \left\{ \begin{array}{ll} \frac{2
      \sqrt{2}(\pi+1) H_3}{h_3 |\zeta_3|}, \quad &\mbox{if}~ \vzeta
    \in \Lambda_\tau^+, \\ \Big[ \delta_{\zeta_1,0} +
      \frac{2H(1-\delta_{\zeta_1,0})}{h |\zeta_1|}\Big] \Big[
      \delta_{\zeta_2,0} + \frac{2H(1-\delta_{\zeta_2,0})}{h
        |\zeta_2|}\Big], & \mbox{if}~ \vzeta \in \Lambda_\tau^0.
  \end{array}
  \right.
  \label{eq:PP9}
\end{equation}
Then we estimate the sum in the right hand side of (\ref{eq:PP8}).  To
prove (\ref{eq:PP8}) we show:
\begin{enumerate}
\item[(i)] {$\Lambda_\tau^+$} contains at least $s$
  grid points.
\item[(ii)] For any $\vzeta \in \Lambda_\tau$, we have
  \[
  \cU\Big(\frac{h \zeta_{1}}{H},\frac{h_3 \zeta_{3}}{H_3}\Big)
  \cU\Big(\frac{h \zeta_{2}}{H},\frac{h_3 \zeta_{3}}{H_3}\Big) <
  \mathcal{Z}(\vzeta).
  \]
\item[(iii)] For any $\vzeta \notin \Lambda_\tau$,
\[
  \cU\Big(\frac{h \zeta_{1}}{H},\frac{h_3 \zeta_{3}}{H_3}\Big)
  \cU\Big(\frac{h \zeta_{2}}{H},\frac{h_3 \zeta_{3}}{H_3}\Big) \le
  \frac{2 \sqrt{2}(\pi+1)}{\tau} = \min_{{\vzeta \in
    \Lambda_\tau^+}} \mathcal{Z}(\vzeta).
  \]
\end{enumerate}
By showing $(i)$, we establish that there are at least as many terms
to sum on the right hand side of (\ref{eq:PP8}) as on the left, and we
can define a one to one map $\mathscr{M}:\Lambda_\star \to
\Lambda_\tau$, such that $\mathscr{M}(\vzeta) = \vzeta$ if $\vzeta \in
\Lambda_\star \cap \Lambda_\tau$ {and
  $\mathscr{M}(\vzeta) \in \Lambda_\tau^+$ if $\vzeta \in
  \Lambda_\star \setminus \Big(\Lambda_\star \cap
  \Lambda_\tau\Big)$}. Points $(ii)$ and $(iii)$ ensure that
\[
\cU\Big(\frac{h \zeta_{1}}{H},\frac{h_3 \zeta_{3}}{H_3}\Big)
\cU\Big(\frac{h \zeta_{2}}{H},\frac{h_3 \zeta_{3}}{H_3}\Big) \le
\mathcal{Z}(\mathscr{M}(\vzeta)), \qquad \forall \vzeta \in
\Lambda_\star,
\]
from which (\ref{eq:PP8}) follows.

\vspace{0.1in} \noindent \textbf{Proof of $(i)$:} This statement is
trivial when $\tau = {s h}/H$ or $\tau = s h_3/H_3$, because
 $\Lambda_\tau^+$ has
cardinality larger than $s$, by definition.  Thus, let $\tau =
\Big[\frac{3s}{8} \Big(\frac{h}{H}\Big)^2 \frac{h_3}{H_3}\Big]^{1/3}
$, and calculate the number $s_\tau$ of grid points in $\Lambda_\tau^+$
as
\begin{align*}
  s_\tau =  \sum_{\vzeta \in
    \Lambda_\tau^+} 1 & \ge \hspace{-0.1in} \sum_{\zeta_3 \in
    \mathbb{Z}, 0 < |\zeta_3| \le \tau H_3/h_3} \left[2 \left \lfloor
    \frac{H}{h} \left( \frac{\tau}{\sqrt{h_3 |\zeta_3|/H_3}} +
    \frac{h_3}{H_3}|\zeta_3|\right) \right \rfloor \right]^2 \\ & \ge 4
    \Big(\frac{H h_3}{h H_3}\Big)^2 \hspace{-0.1in}\sum_{\zeta_3 \in
      \mathbb{Z}, 0 < |\zeta_3| \le \tau H_3/h_3} |\zeta_3|^2 \\ &\ge
    4 \Big(\frac{H h_3}{h H_3}\Big)^2 \frac{2}{3} \Big( \frac{\tau
      H_3}{h_3}\Big)^3 \\ & = \frac{8}{3} \Big(\frac{H}{h}\Big)^2
    \frac{H_3}{h_3} \tau^3 \ge s.
\end{align*}
The first inequality is by definition of $\Lambda_\tau^+$. The
factor $2$ is due to the absolute values and $\lfloor \cdot \rfloor$
denotes the integer part. The second inequality is because we omit one
positive term in the sum. The third inequality is by direct summation
\[
\sum_{j=1}^n j^2 = \frac{n(n+1)(2n+1)}{6} > \frac{n^3}{3}.
\]
Again the factor $2$ is due to the absolute values. The last
inequality is by definition of $\tau$. This concludes the proof of
$(i)$.

\vspace{0.1in}
\noindent \textbf{Proof of $(ii)$:} This follows immediately from
bounds (\ref{eq:PP1}) and (\ref{eq:PP2}).

\vspace{0.1in}
\noindent \textbf{Proof of $(iii)$:} We note from definition
(\ref{eq:PP9}) of $\mathcal{Z}(\vzeta)$ and  (\ref{eq:PP7})
that
\[
\min_{\vzeta \in \Lambda_\tau^+} \mathcal{Z}(\vzeta) = \frac{2
  \sqrt{2}(\pi+1)}{\tau}.
\]
Now consider an arbitrary $\vzeta \notin \Lambda_\tau$. Recalling
definitions (\ref{eq:PP7})-(\ref{eq:PP6}), we distinguish three cases:

\noindent 1. If $|\zeta_3| \ge \frac{\tau H_3}{h_3}$ we have by
(\ref{eq:PP2}) that
  \[
  \cU\Big(\frac{h \zeta_{1}}{H},\frac{h_3 \zeta_{3}}{H_3}\Big)
  \cU\Big(\frac{h \zeta_{2}}{H},\frac{h_3 \zeta_{3}}{H_3}\Big) \le
  \frac{8}{h_3 |\zeta_3|/H_3} \le \frac{8}{\tau} < \frac{2
    \sqrt{2}(\pi+1)}{\tau}.
  \]

\noindent 2.  If $\zeta_3 = 0$ we can assume without loss of
generality that $|\zeta_1| > \frac{\tau}{h/H}$ since at least one of
$\zeta_1$ and $\zeta_2$ must satisfy this condition. We obtain from
(\ref{eq:PP1}) that
  \begin{align*}
  \cU\Big(\frac{h \zeta_{1}}{H},0\Big) \cU\Big(\frac{h
    \zeta_{2}}{H},0\Big) \le \frac{2}{h |\zeta_1|/H} < \frac{2}{\tau}
  < \frac{2 \sqrt{2}(\pi+1)}{\tau}.
  \end{align*}

\noindent 3. If $0< |\zeta_3| \le \frac{\tau H_3}{h_3}$ we can assume
without loss of generality that
  \[
  |\zeta_1| > \frac{H}{h} \left[ \frac{\tau}{\sqrt{h_3 |\zeta_3|/H_3}} +
    \frac{h_3}{H_3} |\zeta_3| \right],
  \]
  since at least one of $\zeta_1$ and $\zeta_2$ must satisfy this
  condition. Then 
  \begin{align*}
  \cU\left( \frac{h}{H} \zeta_1, \frac{h_3}{H_3} \zeta_3 \right) \le
  \frac{\pi+1}{\tau} \sqrt{\frac{h_3}{H_3} |\zeta_3|} ,
  \end{align*}
  by estimate (\ref{eq:PP10}), and using (\ref{eq:PP2}) for the other
  Fresnel integral we get
  \[
  \cU\Big(\frac{h \zeta_{1}}{H},\frac{h_3 \zeta_{3}}{H_3}\Big)
  \cU\Big(\frac{h \zeta_{2}}{H},\frac{h_3 \zeta_{3}}{H_3}\Big) \le
  \frac{2 \sqrt{2} (\pi+1)}{\tau} .
  \]
This concludes the proof of (\ref{eq:PP8}).

It remains to estimate the right hand side in (\ref{eq:PP8}),
\begin{align}
  \sum_{\vzeta \in \Lambda_\tau} \mathcal{Z}(\vzeta) = &
  \frac{2\sqrt{2} (\pi+1)H_3}{h_3} \sum_{\vzeta \in \Lambda_\tau^+}
  \frac{1}{|\zeta_3|} + \sum_{\vzeta \in \Lambda_\tau^0}
  \Big[ \delta_{\zeta_1,0} + \frac{2H(1-\delta_{\zeta_1,0})}{h
      |\zeta_1|}\Big] \Big[ \delta_{\zeta_2,0} +
    \frac{2H(1-\delta_{\zeta_2,0})}{h |\zeta_2|}\Big].
  \label{eq:PP11}
    \end{align}
For the first sum we have by the definition (\ref{eq:PP7}) of
$\Lambda_\tau^+$ that
\begin{align*}
   \frac{H_3}{h_3} \sum_{\vzeta \in \Lambda_\tau^+}
   \frac{1}{|\zeta_3|} &\le \frac{H_3}{h_3}\hspace{-0.03in}
   \sum_{\zeta_3 \in \mathbb{Z}, 0 < |\zeta_3| \le \frac{H_3
       \tau}{h_3}} \frac{1}{|\zeta_3|} \left[1 + \frac{2 H}{h} \left(
     \frac{\tau}{\sqrt{h_3 |\zeta_3|/H_3}} + \frac{h_3}{H_3}
     |\zeta_3|\right) \right]^2 \\  &\le \hspace{-0.1in}\sum_{\zeta_3
     \in \mathbb{Z}, 0 < |\zeta_3| \le \frac{H_3
       \tau}{h_3}} \hspace{-0.05in}\left\{16 \Big(\frac{H}{h}\Big)^2
   \left[ \frac{\tau^2}{\big(h_3 |\zeta_3|/H_3\big)^2} +
     \frac{h_3}{H_3} |\zeta_3| \right] + \frac{2 H_3}{h_3 |\zeta_3|}
   \right\}.
\end{align*}
Now using that
\[
\sum_{j=1}^n j \le \frac{n(n+1)}{2}, \qquad \sum_{j=1}^n \frac{1}{j^2}
\le \frac{\pi^2}{6}, \qquad \sum_{j=1}^n \frac{1}{j} \le 1 + \ln(j),
\]
and substituting in the bound above, we get
\begin{align*}
  \frac{H_3}{h_3} \sum_{\vzeta \in \Lambda_\tau^+} \frac{1}{|\zeta_3|}
  &\le 8 \Big(\frac{H}{h}\Big)^2 \left[ \tau + \frac{\tau^2 H_3}{h_3}
    \Big(1+ \frac{\pi^2 H_3}{3 h_3}\Big)\right] + \frac{2 H_3}{h_3}
  \left[ 1 + \ln\Big(\frac{\tau H_3}{h_3}\Big) \right].
\end{align*}
For the second sum in (\ref{eq:PP11}) we have
\begin{align*}
  \sum_{\vzeta \in \Lambda_\tau^0} \Big[ \delta_{\zeta_1,0} +
    \frac{2H(1-\delta_{\zeta_1,0})}{h |\zeta_1|}\Big] \Big[
    \delta_{\zeta_2,0} + \frac{2H(1-\delta_{\zeta_2,0})}{h
      |\zeta_2|}\Big] = \sum_{\vzeta \in \Lambda_\tau^0} \Big[
    \delta_{\zeta_1,0} (1-\delta_{\zeta_2,0}) \frac{2 H}{h |\zeta_2|}
    + \\\delta_{\zeta_2,0} (1-\delta_{\zeta_1,0}) \frac{2 H}{h
      |\zeta_1|} + (1-\delta_{\zeta_1,0})
    (1-\delta_{\zeta_2,0})\frac{2 H}{h |\zeta_1|} \frac{2H}{h
      |\zeta_2|} \Big],
  \end{align*}
because we cannot have both $\zeta_1$ and $\zeta_2 = 0$. The first
term is bounded as
\[
\sum_{\vzeta \in \Lambda_\tau^0} \delta_{\zeta_1,0}
  (1-\delta_{\zeta_2,0}) \frac{2 H}{h |\zeta_2|} = \frac{4H}{h}
  \sum_{\zeta=1}^{\lfloor \tau H/h\rfloor} \frac{1}{\zeta} \le
  \frac{4H}{h} \big[ 1 + \ln \big(\tau H/h\big)\big],
\]
and similar for the second term. For the last term we have
\begin{align*}
\sum_{\vzeta \in \Lambda_\tau^0} (1-\delta_{\zeta_1,0})
(1-\delta_{\zeta_2,0})\frac{2 H}{h |\zeta_1|} \frac{2H}{h |\zeta_2|} &
=4 \Big( \frac{2 H}{h}\Big)^2 \left(\sum_{\zeta=1}^{\lfloor \tau
  H/h\rfloor} \frac{1}{\zeta} \right)^2 \le \frac{16 H^2}{h^2}
\big[ 1 + \ln \big(\tau H/h\big)\big]^2.
\end{align*}
Gathering the results we have for large $s$, and therefore large
$\tau$,
\begin{equation}
  \sum_{\vzeta \in \Lambda_\tau} \mathcal{Z}(\vzeta) \le \widetilde C
  \Big(\frac{H}{h}\Big)^2 \frac{H_3}{h_3} \tau^2 + \widetilde C_1
  \frac{H_3}{h_3} \ln \tau + \widetilde C_2 \left(\frac{H}{h}\right)^2
  \ln^2 \hspace{-0.05in} \tau,
  \label{eq:PP13}
\end{equation}  
with constant $\widetilde C$ close to $16 \sqrt{2} (\pi+1)$,
$\widetilde C_1$ close to $2$ and $\widetilde C_2$ close to $16$. Here
we used the expectation that $h>h^\star = 4 H$ and $h_3 > h_3^\star =
16 H$. To finish the proof of the Lemma we obtain from definition
(\ref{eq:PP5}) of $\tau$ that
\[
\frac{H_3}{h_3} \ln \tau \le \frac{\ln \Big(s h_3/H_3\Big)}{h_3/H_3}
\le \frac{ \ln s}{h_3/H_3} + e^{-1},
\]
and similarly for $H/h \ln \tau$, where we used that $\ln x/x$ attains
its maximum over the interval $[1,\infty)$ at $x = e$. Moreover, we
  note that the first term in (\ref{eq:PP13}) is negligible in
  comparison with the others unless $\tau = \Big[\frac{3s}{8}
    \Big(\frac{h}{H}\Big)^2 \frac{h_3}{H_3}\Big]^{1/3}$. Substituting
  this expression of $\tau$ in the first term and using that $s$ is
  large, we get Lemma \ref{lem.2} with constant $C$ close to
  $(3/2)^{2/3}(\pi+1)$, $C_1$ close to $\widetilde C_1$ and $C_2$
  close to $\widetilde C_2$.  This concludes the proof of Theorem
  \ref{thm.2}. $\Box$

\subsubsection{Broad band}
\label{sect:parax_pf_bb}

The proof of Theorem \ref{thm.3} is similar to that of Theorems
\ref{thm.1} and \ref{thm.2}, with modifications that account for the
faster decay with the range offset of the inner products
(\ref{eq:BB7}) of the columns of the sensing matrix.

Using the translation invariance of (\ref{eq:BB7}) and writing
explicitly that the points $\vz_q$ are on the grid, we can write the
cumulative coherence as
\begin{equation}
  \mu(\bG,s) = \max_{|\Lambda| = s-1} \sum_{\vzeta \in \Lambda}
  e^{-\big(\frac{h_3 |\zeta_3|}{\cH_3}\big)^2} \Big|
  \mbox{sinc}\Big(\frac{h |\zeta_1|}{\cH} \Big) \Big| \Big|
  \mbox{sinc}\Big(\frac{h |\zeta_2|}{\cH} \Big) \Big|.
  \label{eq:PfBB1}
\end{equation}
Here $\Lambda \in \mathbb{Z}^3 \setminus \{0\}$ is a set of
cardinality $s-1$, as before, and we introduced the notation
\[
\cH_3 =\frac{\sqrt{2}c}{B}, \qquad \cH = \frac{2 L}{k_o a}.
\]

To derive the base resolution limits, let $s = 2$ in (\ref{eq:PfBB1})
and observe that
\begin{equation}
  e^{-\big(\frac{h_3 |\zeta_3|}{\cH_3}\big)^2} \Big|
  \mbox{sinc}\Big(\frac{h |\zeta_1|}{\cH} \Big) \Big| \Big|
  \mbox{sinc}\Big(\frac{h |\zeta_2|}{\cH} \Big) \Big| \le
  e^{-\big(\frac{h_3 |\zeta_3|}{\cH_3}\big)^2}, \qquad \mbox{for} ~ ~
  \zeta_3 \ne 0, 
  \label{eq:PfBB2}
\end{equation}
uniformly in $\zeta_1, \zeta_2 \in \mathbb{Z}$, and when $\zeta_3 =
0$,
\begin{equation}
  \Big| \mbox{sinc}\Big(\frac{h |\zeta_1|}{\cH} \Big) \Big| \Big|
  \mbox{sinc}\Big(\frac{h |\zeta_2|}{\cH} \Big) \Big| \le
   \Big[ \delta_{\zeta_1,0} +
    \frac{1-\delta_{\zeta_1,0}}{h |\zeta_1|/\cH}\Big]\Big[
    \delta_{\zeta_2,0} + \frac{1-\delta_{\zeta_2,0}}{h
      |\zeta_2|/\cH}\Big].
  \label{eq:PfBB3}
\end{equation}
Thus, $\mu(\bG,2) < 1/2$ when $ {\cH}/{h} < {1}/{2}$ and
$e^{-h_3^2/{\cH_3^2}} < {1}/{2}$ or, equivalently, when
\begin{equation}
  h > 2 \cH = h^\star \quad \mbox{and} \quad h_3 > \sqrt{\ln(2)} \cH_3 =
  h_3^\star,
  \label{eq:PfBB4}
\end{equation}
as stated in equation (\ref{eq:BB8}) of Theorem \ref{thm.3}.

To prove the second statement of Theorem \ref{thm.3}, for large $s$,
we use assumption (\ref{eq:BB9}) and the relation (\ref{eq:PfBB4}) to
write
\begin{equation}
  \frac{h}{\cH}, \frac{h_3}{\cH_3} > \beta := \widetilde C \ln s,
  \label{eq:PfBB5}
\end{equation}
for a constant $\widetilde C$ that is slightly larger than $C$.  We
also define $\mathcal{Z}:\mathbb{Z}^3 \to \mathbb{R}$ by
\begin{equation}
  \mathcal{Z}(\vzeta) = e^{-\beta^2 \zeta_3^2} \Big[
    \delta_{\zeta_1,0} + \frac{1-\delta_{\zeta_1,0}}{ \beta
      |\zeta_1|}\Big]\Big[ \delta_{\zeta_2,0} +
    \frac{1-\delta_{\zeta_2,0}}{\beta |\zeta_2|}\Big],
\label{eq:CP1}
\end{equation}
and using (\ref{eq:PfBB5}) in (\ref{eq:PfBB1}) we get
\begin{equation}
  \mu(\bG,s) \le \max_{|\Lambda| = s-1} \sum_{\vzeta \in \Lambda}
  \mathcal{Z}(\vzeta).
  \label{eq:CP2}
\end{equation}

Now let $\Lambda_\star$ denote the maximizing set in
(\ref{eq:CP2}). We do not know it explicitly, but we can define a
one-to-one mapping from $\Lambda_\star$ to another set $\Lambda_\beta$
which allows us to bound $\mu(\bG,s)$.  The construction of the set
\begin{equation}
  \Lambda_\beta = \Big\{\vzeta \in \mathbb{Z}^3\setminus \{0\} \quad
  \mbox{s.t.} \quad \mathcal{Z}(\vzeta) \ge \frac{1}{\beta s} \Big \},
  \label{eq:CP3}
\end{equation}
is motivated by the rapid decay in range of the terms in the sum in
(\ref{eq:PfBB1}). Explicitly, we note that when $\vzeta \in
\Lambda_\beta$ we have $\zeta_3 = 0$ because if this were not true,
definition (\ref{eq:CP1}) would give
\[
\mathcal{Z}(\vzeta) \le e^{-\beta^2 \zeta_3^2} \le e^{-\beta^2} \le
e^{-2 \beta} < e^{-\ln(\beta s)} = \frac{1}{\beta s}.
\]
Here we assumed $\beta > 2$, which is consistent with (\ref{eq:PfBB5})
for large $s$, and since $\beta > \ln \beta$, we also have $ 2 \beta >
\ln(\beta) + \ln s = \ln (\beta s).$

Let $\Lambda_\beta^j$ be the intersection of $\Lambda_\beta$ with the
$\zeta_j$ axis, for $j = 1, 2$. Then, the cardinality
$|\Lambda_\beta|$ of the set $\Lambda_\beta$ satisfies
\[
  |\Lambda_\beta| \ge |\Lambda_\beta^1| + |\Lambda_\beta^2| = 4 s,
\]
because by definition $(\ref{eq:CP3})$, $\zeta_j \in \Lambda_\beta^j$
means that $|\zeta_j| \le s$, for $j = 1, 2$. Thus, there are at least
$4 s$ points in $\Lambda_\beta$, and we can define a one to one
mapping from the maximizing set $\Lambda_\star$ to
$\Lambda_\beta$. Moreover, since for any $\vzeta \notin \Lambda_\beta$
we have $\mathcal{Z}(\vzeta) < 1/(\beta s)$ we conclude from (\ref{eq:CP2}) that
\[
  \mu(\bG,s) \le \sum_{\vzeta \in \Lambda_\beta} \mathcal{Z}(\vzeta). 
\]
To bound the right hand side in this equation, note from definitions
(\ref{eq:CP1}) and (\ref{eq:CP3}) that $\Lambda_\beta$ is contained in
the punctured disk $\mathcal{D}_s$ of radius $s$,
\[
\mathcal{D}_s =\Big \{ \vzeta \in \mathbb{Z}^3\setminus\{0\} \quad
\mbox{s.t.} \quad \zeta_3 = 0, ~ ~ |\zeta_1|, |\zeta_2| \le s \Big\},
\]
and obtain
\begin{align*}
  \mu(\bG,s) \le \sum_{\vzeta \in \mathcal{D}_s} \mathcal{Z}(\vzeta)
  &= \prod_{j=1}^2 \sum_{|\zeta_j| \le s} \left[ \delta_{\zeta_j,0} +
    \frac{1-\delta_{\zeta_j,0}}{\beta |\zeta_j|} \right] -
  \mathcal{Z}(0) \le \left[ 1 + \frac{2}{\beta} \Big( 1 + \ln s
    \Big) \right]^2 - 1 \\ & = \frac{4 \ln s}{\beta} + \frac{4
    \ln^2 \hspace{-0.05in} s}{\beta^2}  = \frac{4 \ln s}{\beta} +
  \frac{4}{\beta} + \left[\frac{2}{\beta}(1 + \ln s) \right]^2.
\end{align*}
The proof of Theorem \ref{thm.3} is completed with the observation
that we can make the bound in this estimate less than $1/2$ by
choosing the constant $\widetilde C$ in (\ref{eq:PfBB5}) large enough,
independent of $s$. $\Box$


\section{Imaging on fine grids}
\label{sect:cont}
The resolution estimates in Theorems \ref{thm.1}-\ref{thm.3} do not
account for noise and modeling errors due to off-grid placement of the
sources, which may be large for coarser discretizations required by
the theorems.
\begin{figure}[t]
\begin{center}
\includegraphics[width =
  0.33\textwidth]{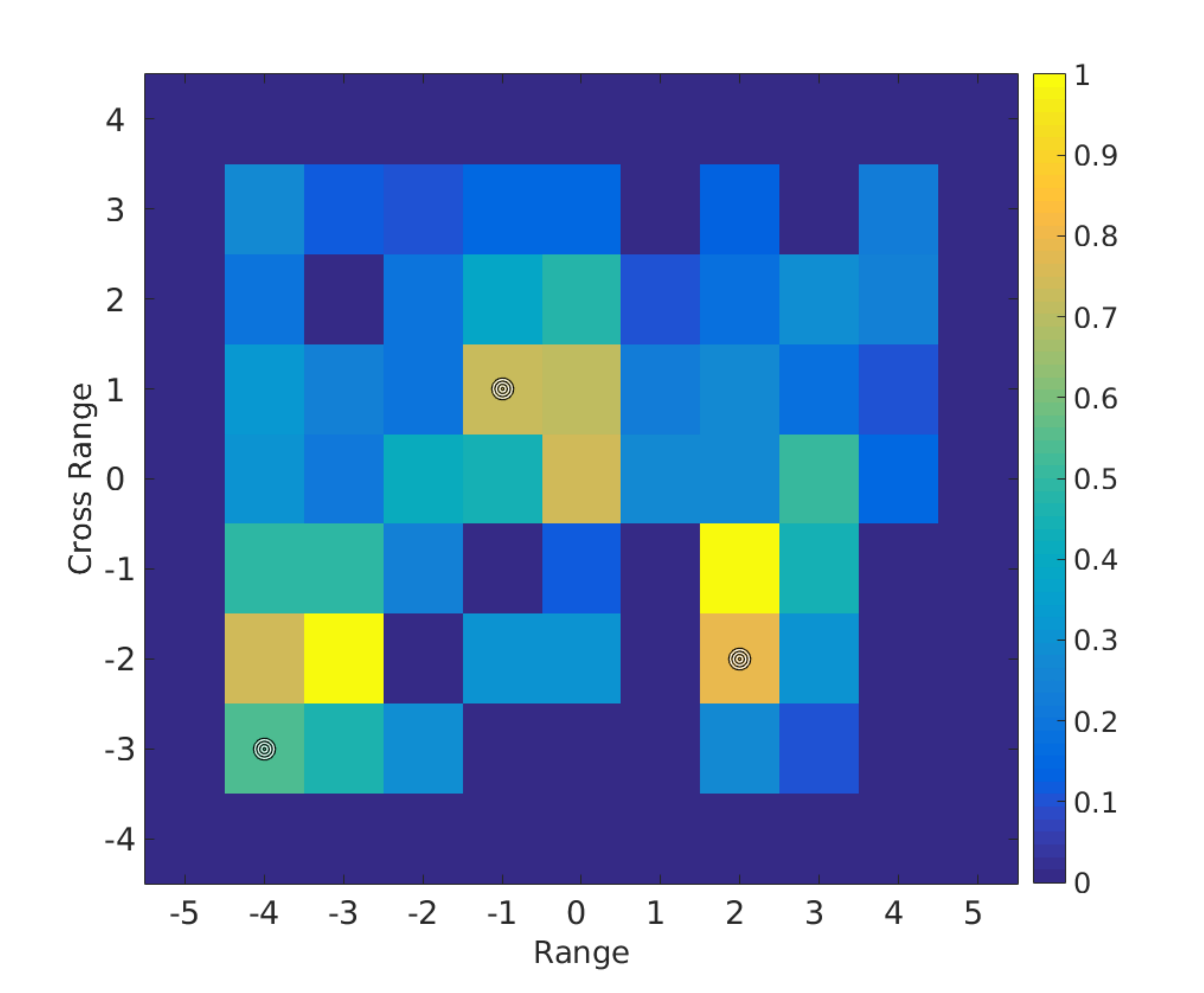}
\hspace{-0.2in}\includegraphics[width =
  0.33\textwidth]{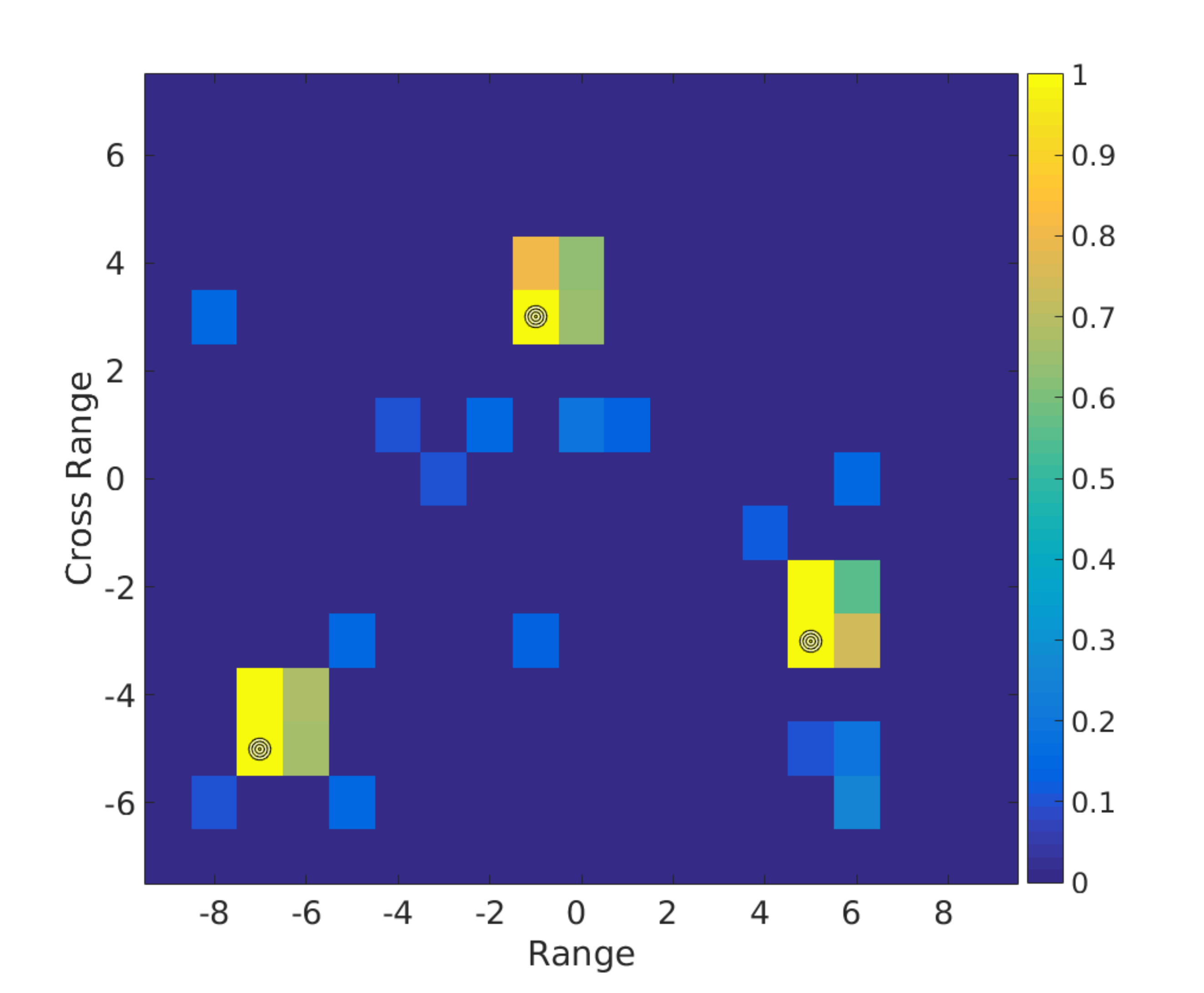}
\hspace{-0.2in}\includegraphics[width =
  0.33\textwidth]{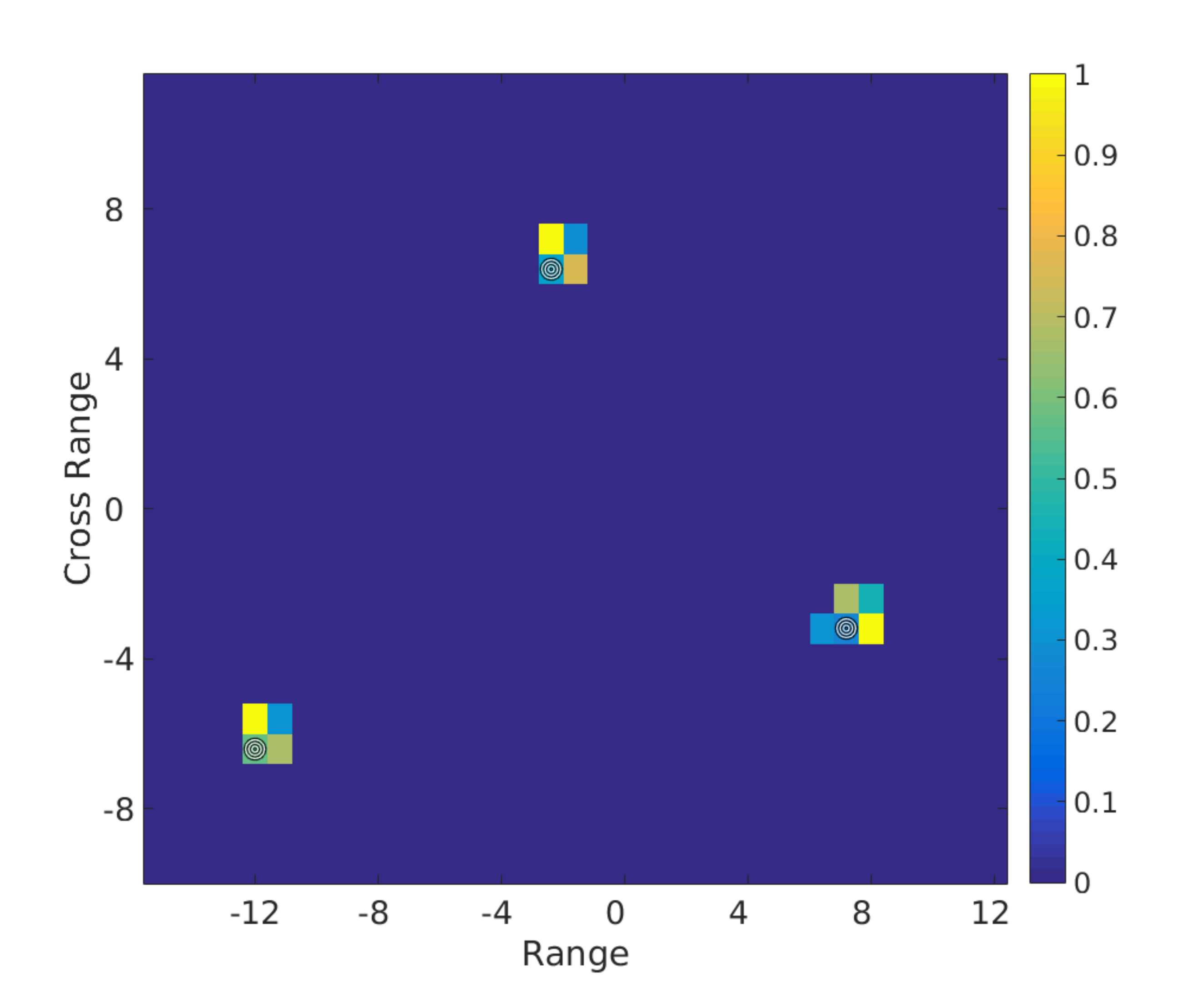}
\end{center}

\caption{Effect of modeling error due to sources off-grid. From left
  to right: Image discretized at base resolution, at $1/2$ base
  resolution and $1/4$ base resolution.}
\label{fig:trade-off}
\end{figure}
In this section we mitigate the modeling error by discretizing the
imaging region on a fine mesh, and then study the results of the
$\ell_1$ optimization. As the results in the previous section shows it
is impossible to have a meaningful answer for arbitrary distributions
of the sources. However, if they are located at points or clusters of
points that are sufficiently far apart, the results are useful, as
illustrated in Figure \ref{fig:trade-off}. In the left image we
display the result of the optimization for a discretizations of the
imaging region at the base resolution $\vbh^\star$ defined in Theorem
\ref{thm.1}. If the sources were on the grid, the reconstruction would
have been perfect. Here the sources are off-grid, and the
reconstruction is poor due to the modeling error. The result is better
in the other two plots because the modeling error is reduced by taking
a smaller mesh size.  While in general a discretization at $1/4 \,
\vbh^\star$ does not guarantee a good recovery, here the result is
good because the sources are far apart, and thus have little
interaction. The analysis in this section formalizes this observation.
We point the reader to \cite{fannjiang2012coherence} for a different
study of similar ideas, and algorithms designed to take advantage of
the weak interaction between the sources.  Here we study directly the
optimization problems (\ref{eq:F7}) and (\ref{eq:F7b}) and consider in
addition clusters of sources.

We begin in section \ref{sect:FGStat} with the statement of results
for well separated sources and then consider clusters of sources in
section \ref{sect:FGStatClust}. The proofs are in section
\ref{sect:FGProof}.

\subsection{Statement of results for well separated sources}
\label{sect:FGStat}
Let us modify the notation slightly, and call $\bg_{\vy}$ the
normalized vector of Green's functions taking us from the array to the
point $\vy$ in the imaging region. When $\vy$ is a point $\vz_j$ on
the grid, then $\bg_{\vy}$ is the same as $\bg_j$ defined before, and
we have the simpler notation $\bg_j \equiv \bg_{\vz_j}$. We quantify
the interaction between two sources located at $\vy$ and $\vy'$ in the
imaging region $W$ by the value of $\Big| \lin \bg_{\vy}, \bg_{\vy'}
\rin \Big|$. Explicitly, in terms of the semi-metric $\D:W\times W \to
     [0,1]$,
\begin{equation}
  \D(\vy,\vy') = 1- \Big| \lin \bg_{\vy}, \bg_{\vy'} \rin \Big|, \qquad
  \forall \, \vy, \vy' \in W,
  \label{eq:FG1}
\end{equation}
that defines the open ball 
\begin{equation}
  \cB_r(\vy) = \Big\{ \vy' \in \mathbb{R}^3 \quad \mbox{s.t.} \quad
  \D(\vy,\vy') < r \Big\},
  \label{eq:FG2}
\end{equation}
we say that points $\vy'$ outside $\cB_r(\vy)$ have a weaker
interaction with $\vy$ than the points in the ball,
\begin{equation}
  \Big| \lin \bg_{\vy}, \bg_{\vy'} \rin \Big| \le 1 - r, \qquad
  \forall \, \vy' \notin \cB_r(\vy).
  \label{eq:FG3}
\end{equation}
In these definitions it does not matter if we have single or multiple
frequency measurements\footnote{The theory applies to both single
  frequency and multiple frequency measurements, but the numerical
  simulations are for a single frequency.}. In both cases we know from
the previous section that $\Big| \lin \bg_{\vy}, \bg_{\vy'} \rin
\Big|$ is a function of $\vy-\vy'$ which peaks at the origin, and is
monotonically decreasing in its vicinity. This means that there exists
a small enough $\bar r$ such that if $r < \bar r$ and $\vy' \in
\cB_r(\vy)$, $\vy$ and $\vy'$ are close in Euclidian distance.

Suppose that we have $s$ sources in the imaging region $W$, supported
at points in the set $\cY = \{\vy_j, ~ ~ j = 1, \ldots s\},$ and
discretize $W$ on a grid with $N$ points $\vz_j$ and mesh size $\vbh$
that is as small as needed to mitigate the modeling error. We define
the interaction coefficient of the set $\cY$ by
\begin{equation}
  \cI(\cY) = \max_{q = 1, \ldots, N} \sum_{\vy_j \in \cY
    \setminus \{\cN(\vz_q)\}} \Big| \lin
  \bg_{\vy_j},\bg_{q} \rin \Big|,
  \label{FG4}
\end{equation}
where $\cN(\vz_q)$ is the closest point to $\vz_q$ in $\cY$, with
respect to semi-metric $\D$. 
{ Here it is possible to make  $\cI(\cY)$ independent of the mesh by replacing the maximum with supremum over the imaging window which is an open subset of $\mathbb{R}^3$.}
Note that $\cI(\cY)$ is similar to the
cumulative mutual coherence $\mu(\bG,s)$, except that the set $\cY$
is fixed and the points in it may not be on the grid.  Note also that
$\cI(\cY)$ is well defined even when there are multiple points in
$\cY$ that are closest to $\vz_q$. In such cases we let $\cN(\vz_q)$
be any one of these points without affecting the value of
$\cI(\cY)$.

The next two theorems describe the support of the $\ell_1$
minimizer. Theorem \ref{thm.4} and its corollary are for formulation
(\ref{eq:F7}) of the optimization problem, which assumes an exact
model. Theorem \ref{thm.5} is for formulation (\ref{eq:F7b}) which
accounts for noise and modeling error.

\vspace{0.05in}
\begin{theorem}
\label{thm.4}
Suppose that the unknown sources are supported on the fine grid at
points $\vz_j$ enumerated by the set $S$ of cardinality $s$, and are
represented by the $s$-sparse vector $\brho\in \complex^N$, satisfying
$\bG \brho = \bd$. The sources are assumed sufficiently far apart so
that for some $r \in (0,1)$ the balls $\cB_r(\vz_j)$ are
disjoint. Let $\brho_\star$ be the solution of the optimization
problem (\ref{eq:F7}) and decompose it as
\begin{equation}
    \brho_\star = \brho_\star^{(i)} + \brho_\star^{(o)},
    \label{eq:FG5}
\end{equation}
where $ \supp \brho_\star^{(i)} \subset \displaystyle
\bigcup_{j\in S} \cB_r(\vz_j), $ and $\brho_\star^{(o)}$ is supported
in the complement of this union. Then,
\begin{equation}
  \|\brho_\star^{(o)}\|_{1} \le \frac{2 \cI(\cY)}{r} \|\brho_\star\|_1.
  \label{eq:FG6}
\end{equation}
\end{theorem}

\vspace{0.05in}
This theorem says that when the interaction coefficient $\cI(\cY)$
is smaller than $r/2$, the support of the optimal $\brho_\star$ is
concentrated in the vicinity of the sources. The next corollary
quantifies the error of the reconstruction.

\vspace{0.05in}
\begin{corollary}
  \label{cor.thm4}
  Under the same assumptions as in Theorem \ref{thm.4}, the error of
  the reconstruction is quantified by
  \begin{equation}
    \|\brho - \bar \brho_\star \|_1 \le \frac{2
        \cI(\cY)}{r} \|\brho\|_1, 
  \label{eq:FG7}
  \end{equation}
  where $\bar \brho_\star$ is the effective source vector in
  $\complex^N$ with $j-$th component given by
  \begin{equation}
    \bar \rho_{\star j} = \left\{ \begin{array}{ll} \displaystyle
      \sum_{q \in \mathscr{S}_j} \rho_{\star q}^{(i)} \lin \bg_j,
      \bg_{q} \rin, \qquad & \mbox{for} \, j \in \cS, \\ 0, & \mbox{for}
      \, j \notin \cS,
      \end{array} \right.
    \label{eq:FG8}
  \end{equation}
  where $\rho_{\star q}^{(i)}$ denotes the $q-$th component of
  $\brho_\star^{(i)}$ and $\mathscr{S}_j$ is the set\footnote{Note
    that the set $\mathscr{S}_j$ depends on $r$, but for simplicity we
    suppress the dependence in the notation.} of indexes of the grid
  points supported in $\cB_r(\vz_j)$.
\end{corollary}

\vspace{0.05in}
Note that $\bar \brho_\star$ is an $s$-sparse vector of the same
support $\cS$ as $\brho$, but with entries given by the ``weighted''
sum of the components of $\brho_\star^{(i)}$ supported in the vicinity
of each source. When the radius $r$ is small, the complex weights
$\lin \bg_j, \bg_q \rin$ are close to one, and $\bar \rho_{\star j}$
is approximately the sum of the components of $\brho_\star^{(i)}$
supported in $\mathscr{S}_j$.
Furthermore, (\ref{eq:FG7}) implies that when $ \frac{2
        \cI(\cY)}{r}$ is sufficiently small such that the right hand side of (\ref{eq:FG7}) is less than  $\min_{j \in \cS} \{|\rho_j|\} $ ,   $\brho_\star$ has a non-zero component in the $r$-neighborhood of every source location.   
\begin{figure}[t]
  \begin{center}
    \vspace{-0.05in}
    \hspace{-0.05in}\includegraphics[width=0.5\textwidth]{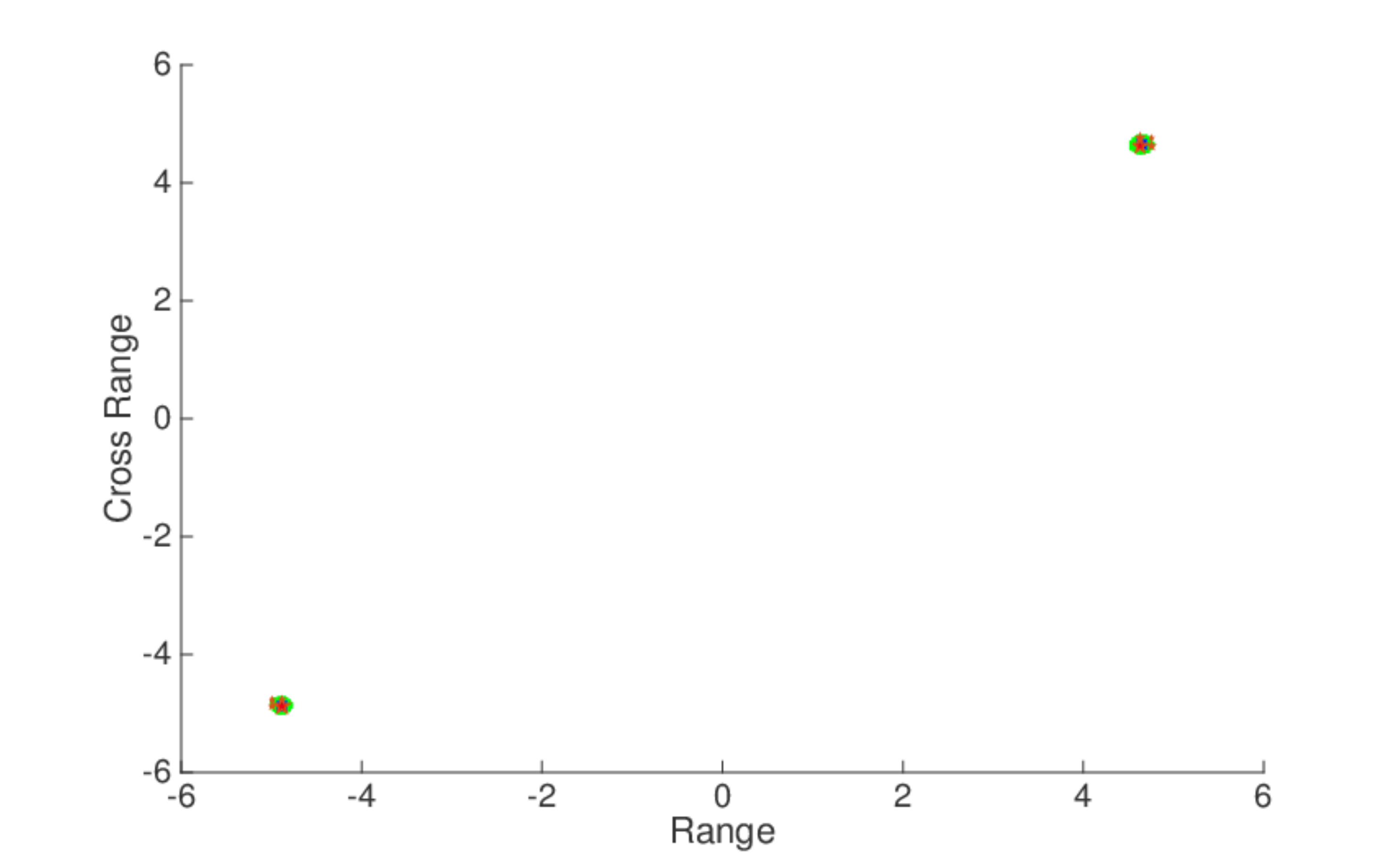}
\includegraphics[width=0.5\textwidth]{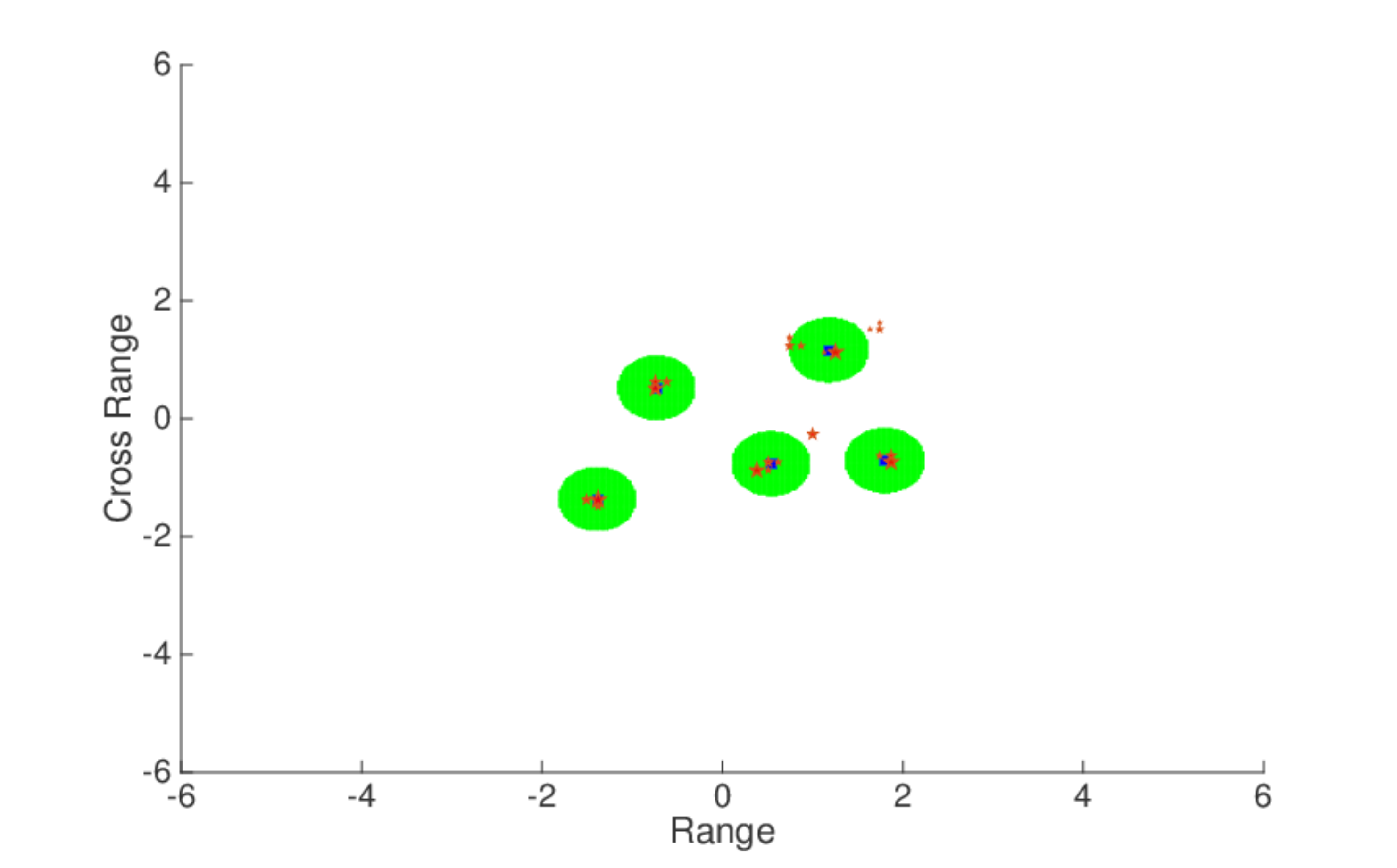}
  \end{center}

  \vspace{-0.1in}
  \caption{Reconstructions of two sources (left) and five sources
    (right).  The support of $\brho_\star$ is indicated with a star
    of size proportional to its magnitude. The balls $\cB_r(\vy_j)$
    are drawn in green. The radius is $0.009$ in the left plot and
    $0.11$ in the right plot. The axes are range and cross-range in
    units of base resolution $h_3^\star$ and $h^\star$.}
  \label{fig:wellSep1}
\end{figure}

\begin{figure}[t]
  \begin{center}
    \includegraphics[width=0.5\textwidth]{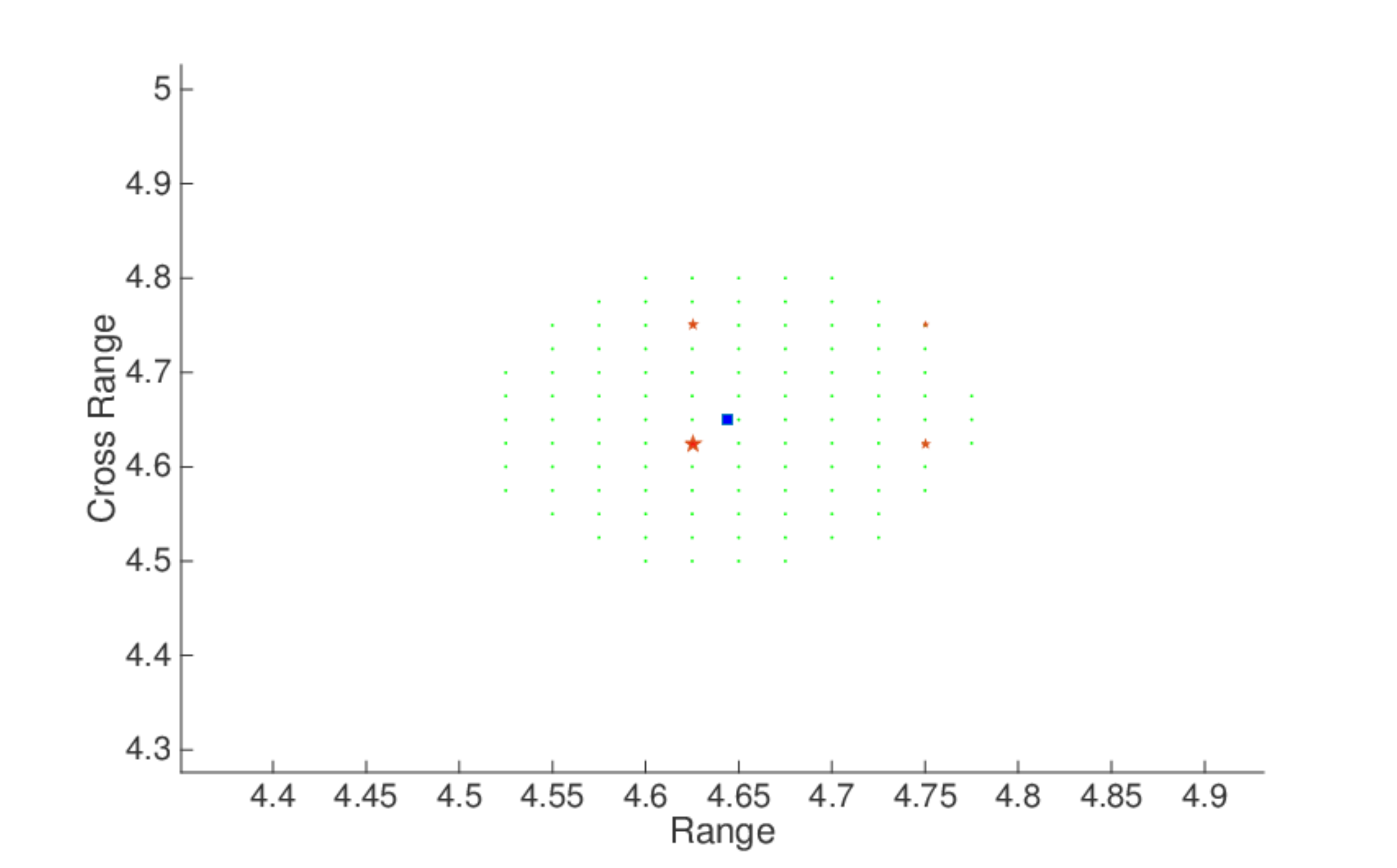}
        \end{center}

  \vspace{-0.1in}
  \caption{Zoom of the image displayed in the left plot of Figure
    \ref{fig:wellSep1} around one of the sources shown with a blue
    square. }
        \label{fig:wellSep2}
\end{figure}

The statements of Theorem \ref{thm.4} and Corollary \ref{cor.thm4} are
already illustrated in the right plot of Figure \ref{fig:trade-off}.
Additional examples are in Figure \ref{fig:wellSep1}, where we show
numerical reconstructions for two sources (left plot) and five sources
(right plot). We display the balls $\cB_r(\vy_j)$ in green, and the
entries in $\brho_\star$ with stars of size proportional to their
magnitude.  In the left plot (see also the zoom in Figure
\ref{fig:wellSep2}) the sources have weak interaction $\cI(\cY) =
0.086$, and for $r = 0.009$ we have
$\|\brho_\star^{(o)}\|_1/\|\brho_\star\|_1 = 2.4\%$ and $\|\brho-\bar
\brho_\star\|_1 / \|\brho\|_1 = 18\%$. For $r = 0.005$ the error drops
to $1.14\%$, however $\|\brho_\star^{(o)}\|_1/\|\brho_\star\|_1$ grows
to $19\%$. That is to say, roughly $80\%$ of the amplitude of the
reconstruction $\brho_\star$ is accumulated very close near the
sources.  In the right plot the interaction coefficient is larger
$\cI(\cY) = 1.43$, but the reconstruction is still good,
$\|\brho_\star^{(o)}\|_1/\|\brho_\star\|_1 = 0.5\%$ and $\|\brho-\bar
\brho_\star\|_1 / \|\brho\|_1 = 10\%$ for $r = 0.11$.

Note that in both simulations the support of the reconstruction is
much better than predicted by Theorem \ref{thm.4}, which gives a
pessimistic bound for $r > 2 \cI(\cY)$. A sharper estimate may be
obtained under the additional assumption that all the entries in
$\brho$ are positive, by taking advantage of cancellations of the
oscillatory terms in the sums analyzed in the proof of the theorem in
section \ref{sect:FGProof1}. However, this is difficult to do without
making strong assumptions on the geometric distribution of the sources
in the imaging region.

The next theorem considers the more general case of an inexact model,
due to noisy data and sources off the grid, and uses the $\ell_1$
penalty formulation (\ref{eq:F7b}). The result is stronger than in
Theorem \ref{thm.1}, as it states that the minimizer $\brho_\star$ is
exactly supported in the vicinity of the sources, for large enough
penalty parameter $\gamma$. This is somewhat expected, as increasing
$\gamma$ in (\ref{eq:F7b}) means putting more emphasis on having a
smaller $\ell_1$ norm i.e., a sparser solution. What is interesting is
that the support of this sparse solution is guaranteed to be near that
of the unknown sources.  However, increasing $\gamma$ comes at the
cost of a larger residual $\|\bG \brho_\star - \bd \|_2$, and there is
no guarantee that the error of recovery of $\brho$ is small, as in
Corollary \ref{cor.thm4}.

\vspace{0.05in}
\begin{theorem}
  \label{thm.5}
  Consider $s$ sources supported in the set $\cY = \{\vy_j, ~ j = 1,
  \ldots, s\}$, with interaction coefficient $\cI(\cY)< 1/2$, so that
  we can find an $r \in (0,1)$ satisfying $r > 2 \cI(\cY)$. Then, for
  sufficiently large penalty parameter $\gamma$, that depends on the
  noise and modeling error, the minimizer $\brho_\star$ of
  (\ref{eq:F7b}) is supported in $\bigcup_{j = 1}^s \cB_r(\vy_j)$.
\end{theorem}

\subsection{Statement of results for clusters of sources}
\label{sect:FGStatClust}
Here we give the generalization of Theorem \ref{thm.4} to clusters of
sources. We assume for simplicity, as in Theorem \ref{thm.4}, that the
sources are on the grid. The result extends to sources off the grid by
modifying the proof of Theorem \ref{thm.5}.

Let us define the effective support $\cS_\ep\subset \cS$ of $\brho \in
\complex^N$, for some $\ep \in (0,1)$, so that
\begin{equation}
  \cY = \{ \vy_1, \ldots, \vy_s \} = \{ \vz_j, ~ j \in \cS \}
  \subset \bigcup_{q \in \cS_\ep } \cB_\ep(\vz_q),
\label{eq:CL1}
\end{equation}
where
\[
  \cB_\ep(\vz_q) \cap \cB_\ep(\vz_l) = \emptyset ~ ~ \forall \, l,
  q \in \cS_\ep, l \ne q.
\]
More explicitly, we cover the set $\cY$ of locations of the sources
with disjoint balls of radius $\ep$, centered at points in $\cS_\ep$.
The set $\cS_\ep$ is the support of the effective source vector
$\bar \brho$, with entries defined similarly to (\ref{eq:FG8})
\begin{equation}
  \bar \rho_{j} = \left\{ \begin{array}{ll} \displaystyle \sum_{q
      \in S \cap \cB_\ep(\vz_j)} \rho_q \lin \bg_q,\bg_j \rin, \quad &
    \mbox{for} \, \, j \in \cS_\ep, \\ 0, & \mbox{otherwise},
    \end{array} \right.
    \label{eq:CL2}
\end{equation}
Obviously $\bar \brho$ depends on $\ep$, but we suppress the
dependence in the notation.  When $\ep \ll 1$, meaning that the
sources are tightly clustered around the points in $\cS_\ep$, the
effective source is approximately the sum of the entries of $\brho$
supported in the cluster. When $\ep$ is larger the complex weights in
(\ref{eq:CL2}) can be far from one and oscillatory, so there may be a
lot of cancellations in the sum in (\ref{eq:CL2}).  Cancellations
(destructive interference of sources) can arise for tight clusters as
well, when the entries in $\brho$ in a cluster have opposite signs.

The result stated in the next theorem says that if the clusters are
far apart and there is little destructive interference of the sources
in the cluster, the support of the optimizer $\brho_\star$ of
(\ref{eq:F7}) is concentrated near the sources.

\vspace{0.05in}
\begin{theorem}
  \label{thm.6}
  Suppose that the unknown sources are supported on the grid at points
  enumerated by the set $\cS$, and that there is an $\ep \in
  (0,1)$ for which we can define the effective support $\cS_\ep$. Let
  $\brho_\star$ be the $\ell_1$ minimizer of (\ref{eq:F7}) and
  decompose it as $\brho = \brho_\star^{(i)} + \brho_\star^{(o)}$,
  where $\brho_\star^{(i)}$ is supported in the disjoint union
  $\bigcup_{j \in \cS_\ep} \cB_r(\vz_j)$, for $r$ satisfying $\ep < r <
  1$, and $\brho_\star^{(o)}$ is supported in the complement of this
  union. We have
  \begin{equation}
    \|\brho_\star^{(o)}\|_1 \le \frac{2 \cI(\cY_\ep)}{r}
    \|\brho_\star\|_1 + \frac{\|\brho\|_1 - \|\bar \brho \|_1}{r},
      \label{eq:CL3}
  \end{equation}
  where $\cY_\ep = \{\vz_j, ~ j \in \cS_\ep \}$ is assumed to satisfy
  $\cI(\cY_\ep) < 1$.
\end{theorem}

\begin{figure}[t]
  \begin{center}
    \vspace{-0.8in}   
    \includegraphics[width=0.5\textwidth]{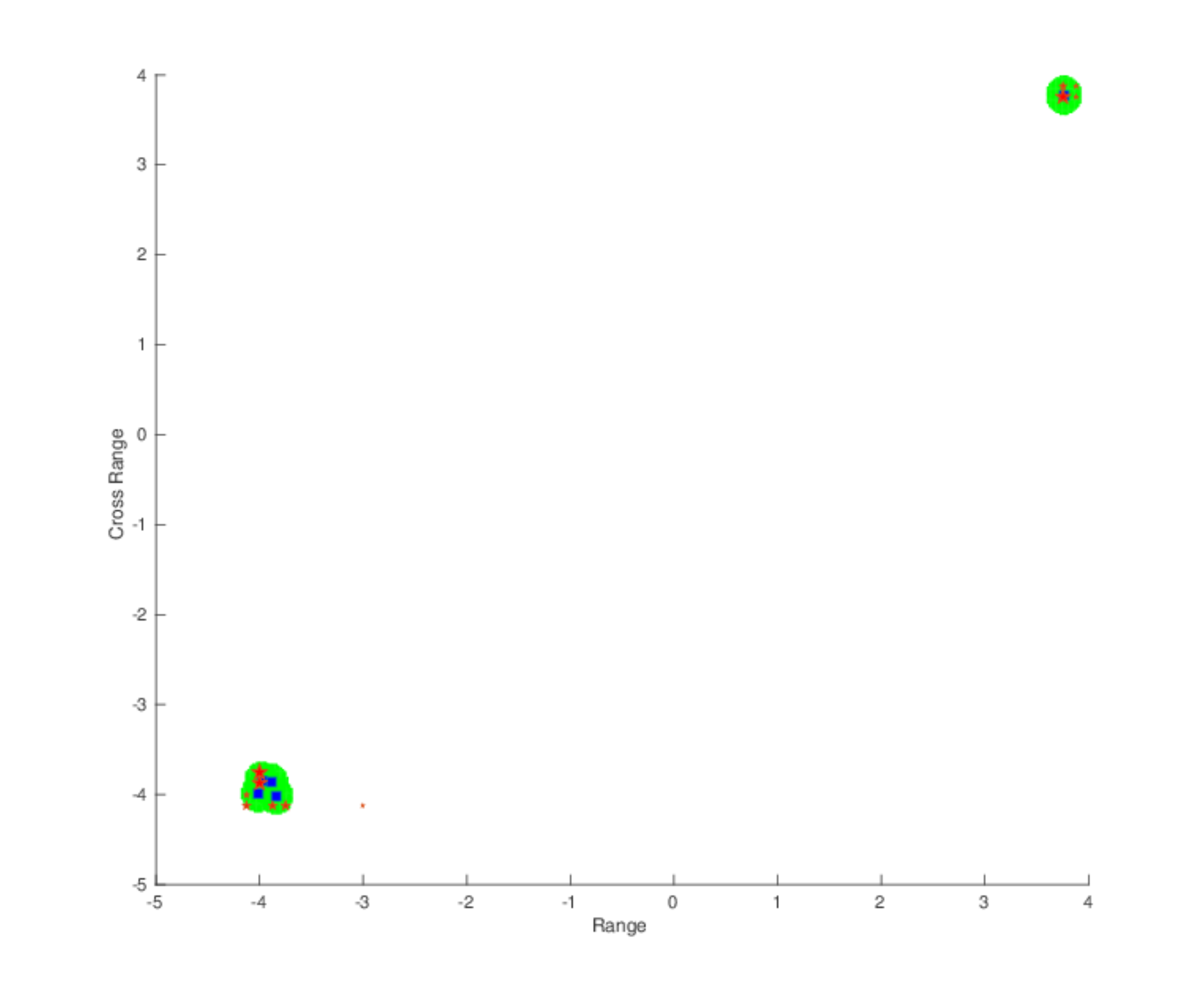}
  \end{center}

\vspace{-0.3in}
  
\caption{Numerical simulation for two clusters of sources. The
    support of $\brho_\star$ is indicated with a star of size
    proportional to its magnitude. The balls $\cB_r$ of radius $r =
    0.017$ are drawn in green. The axes are range and cross-range in
    units of base resolution $h_3^\star$ and
    $h^\star$.  }
  \label{fig:clusters}
\end{figure}

\vspace{0.05in}
The assumption $\cI(\cY_\ep)<1$ is used in the proof, but for the
estimate (\ref{eq:CL3}) to be useful we need $\cI(\cY_\ep) < r/2 <
1/2$. This is because by definition (\ref{eq:CL2}) of $\bar \rho$ we
have $\|\bar \brho\|_1 \le \|\brho \|_1$ and the bound is larger than
$\|\brho_\star\|_1$ for $\cI(\cY_\ep) > r/2$. As before, the estimate
in the theorem is pessimistic.  The numerical results are better, as
illustrated by the simulation in Figure \ref{fig:clusters}. There are
five sources with amplitude equal to one, and locations indicated in
the plot with blue squares. One source is isolated and the other four
form a cluster, so $\cS_\ep$ has cardinality two. The interaction
coefficient of the set $\cY$ is large, $\cI(\cY) = 3.02$, because of
the cluster, but $\cI(\cY_\ep) = 0.096$. The balls shown with
green in the figure are for $r = 0.017$, and the entries in
$\brho_\star$ are indicated with stars of size proportional to the
magnitude. The error is $ \|\brho_\star^{(o)}\|_1/ \|\brho_\star\|_1 =
0.54\% $.  The restriction of $\brho_\star$ to the ball containing the
cluster has $l_1$ norm $2.3$. The restriction to the other ball has
norm $0.96$, which is close to the amplitude of the isolated source.

\subsection{Proofs}
\label{sect:FGProof}
We begin in section \ref{sect:FGProof1} with the proofs of Theorem
\ref{thm.4} and its Corollary \ref{cor.thm4}. Theorem \ref{thm.6} for
clusters of sources is proved in section \ref{sect:FGProof3}. The
proofs of Theorems \ref{thm.4} and \ref{thm.6} are similar but the
proof of Theorem \ref{thm.5} for $\ell_1$ penalty reconstructions is
more involved. We present it in section \ref{sect:FGProof2}.

\subsubsection{$\ell_1$ optimal reconstructions of well separated sources}
\label{sect:FGProof1}
Theorem \ref{thm.4} and its corollary assume an exact model, with well
separated sources on the grid, at points indexed by $\cS$. Recall that
$\mathscr{S}_q$ is the set that enumerates the grid points supported
in $\cB_r(\vz_q)$, for $q \in \cS$. The balls $\cB_r(\vz_q)$ are
disjoint by assumption, so each nonzero entry in $\brho_\star^{(i)}$
is contained in exactly one ball and
\[
\supp \brho_\star^{(i)} = \bigcup_{q \in \cS} \mathscr{S}_q.
\]
By definition of $\brho_\star$ we have $\bG \brho = \bG \brho_\star$
or more explicitly,
\begin{align}
  \sum_{q \in \cS} \rho_q\bg_q = \sum_{q\in \cS} \sum_{j \in
    \mathscr{S}_q} \rho_{\star j}^{(i)} \bg_j + \sum_{j \in
    \mathscr{S}^c} \rho_{\star j}^{(o)} \bg_j,
  \label{eq:FGPf1}
  \end{align}
where we denote the support of $\brho_\star^{(o)}$ by $ \mathscr{S}^c
= \{1, \ldots, N\} \setminus \displaystyle \bigcup_{j \in \cS}
\mathscr{S}_j.  $ The proof of the theorem and its corollary amounts
to estimating the inner products of the left and right sides of
equation (\ref{eq:FGPf1}) with a carefully chosen vector $\bu$, as 
shown next.

\textbf{Proof of Theorem \ref{thm.4}:} Let us define the vector
\begin{equation}
  \bu = \sum_{q\in \cS} \mbox{sign}(\rho_q) \bg_q,
  \label{eq:FGPf2}
\end{equation}
where ``sign'' denotes the complex sign function, and take the inner
product of $\bu$ with the left and right hand side in
(\ref{eq:FGPf1}). We obtain 
\begin{equation}
   \cT_L := \Big| \sum_{q \in \cS} \rho_q \lin \bg_q,\bu \rin \Big| =
   \Big| \sum_{q\in \cS} \sum_{j \in \mathscr{S}_q} \rho_{\star
     j}^{(i)} \lin \bg_j, \bu \rin + \sum_{j \in \mathscr{S}^c}
   \rho_{\star j}^{(o)} \lin \bg_j, \bu \rin \Big| =: \cT_R,
   \label{eq:FGPf3}
\end{equation}
where obviously the left hand side $\cT_L$ equals the right hand side
$\cT_R$. We distinguish them here so we can bound them separately
below and above.

For $\cT_L$ we have
\begin{align*}
  \cT_L & = \Big| \sum_{q \in \cS} \Big[ \rho_q \, \mbox{sign}(\rho_q) +
    \rho_q \sum_{j \in \cS, j \ne q} \mbox{sign}(\rho_j) \lin
    \bg_q,\bg_j\rin \Big| \\ &= \Big| \sum_{q \in \cS} \Big[ |\rho_q
      |+ \rho_q \sum_{j \in \cS, j \ne q} \mbox{sign}(\rho_j) \lin
      \bg_q,\bg_j\rin \Big| \\ & \ge \sum_{q \in \cS} |\rho_q| -
      \sum_{q \in \cS} |\rho_q| \sum_{j \in \cS, j \ne q} \left| \lin
      \bg_q,\bg_j \rin \right|,
 \end{align*}
where we used the definition of the sign function and the triangle
inequality. Since
\[
\sum_{j \in \cS, j \ne q} \left| \lin \bg_q,\bg_j \rin \right| \le
\max_{q \in \cS} \sum_{j \in \cS, j \ne q} \left| \lin \bg_q,\bg_j \rin
\right| \le \max_{q = 1, \ldots, N} \sum_{j \in \cS, j \ne q} \left|
\lin \bg_q,\bg_j \rin \right| = \cI(\cY),
\]
and $\brho$ is supported on $\cS$, we get
\begin{equation}
  \cT_L \ge \|\brho\|_1 \Big[ 1 - \cI(\cY) \Big].
  \label{eq:FGPf4}
\end{equation}

For $\cT_R$ we have by definition (\ref{eq:FGPf2}) of $\bu$ that
\begin{align*}
  \cT_R = \Big| \sum_{q \in \cS} \sum_{j \in \mathscr{S}_q} \rho_{\star
    j}^{(i)} \Big[\mbox{sign}(\rho_q) \lin \bg_j,\bg_q \rin + \sum_{l
      \in \cS, l \ne q}\, \mbox{sign}(\rho_l) \lin \bg_j,\bg_l \rin
    \Big]+   \sum_{j \in \mathscr{S}^c} \sum_{l\in
    \cS} \rho_{\star j}^{(o)} \mbox{sign}(\rho_l) \lin \bg_j, \bg_l \rin
  \Big|
\end{align*}
and using the triangle inequality
\begin{align}
  \cT_R \le \sum_{q \in \cS} \sum_{j \in \mathscr{S}_q} |\rho_{\star
    j}^{(i)}| \left|\lin \bg_j,\bg_q \rin\right| + \sum_{q \in \cS}
  \sum_{j \in \mathscr{S}_q} |\rho_{\star j}^{(i)}| \sum_{l \in \cS, l
    \ne q}\, \left|\lin \bg_j,\bg_l \rin \right|+ \sum_{j \in
    \mathscr{S}^c} |\rho_{\star j}^{(o)}| \sum_{l \in \cS}| \lin
  \bg_j, \bg_l \rin | .
  \label{eq:FGPf5}
\end{align}
In the first term in (\ref{eq:FGPf5}) we have
\[
1 - r < |\lin \bg_j,\bg_q \rin | \le 1,
\]
because $j \in \mathscr{S}_q$. In the second term
\[
\sum_{l \in \cS, l
  \ne q}\, \left|\lin \bg_j,\bg_l \rin \right| \le \cI(\cY),
\]
by definition (\ref{FG4}) of the interaction coefficient and the fact
that $\vz_q$ is the closest source point to $\vz_j$. To bound the
third term in (\ref{eq:FGPf5}), recall that $\cN(\vz_j)$ is
the closest source point to $\vz_j$, for $j \in \mathscr{S}^c$. Its
distance from $\vz_j$ satisfies $\D(\vz_j,\cN(\vz_j)) \ge d$,
by definition of $\mathscr{S}^c$, and therefore
\[
|\lin \bg_j, \bg_{\cN(\vz_j)} \rin| \le 1 - r. 
\]
Moreover,
\[
\sum_{l \in \cS, \vz_l \ne \cN(\vz_j)} |\lin \bg_j, \bg_l \rin|
\le \cI(\cY),
\]
so in the third sum in (\ref{eq:FGPf5}) we have 
\[
\sum_{l \in \cS}| \lin \bg_j, \bg_l \rin | \le \|\brho_\star^{(o)}\|_1
\Big[1-r + \cI(\cY)\Big], \qquad \forall \, j \in \mathscr{S}^c.
\]
Thus, the bound on $\cT_R$ becomes 
\begin{equation}
  \cT_R \le \|\brho_\star^{(i)}\|_1 \Big[1 + \cI(\cY)\Big] +
  \|\brho_\star^{(o)}\|_1 \Big[1-r + \cI(\cY)\Big].
  \label{eq:FGPf6}
\end{equation}

To complete the proof use that $\|\brho_\star\|_1 =
\|\brho_\star^{(i)}\|_1 + \|\brho_\star^{(o)}\|_1$ in (\ref{eq:FGPf6})
and obtain from equations (\ref{eq:FGPf3}) and (\ref{eq:FGPf4}) that
\[
\|\brho_\star\|_1 \Big[1-\cI(\cY)\Big] \le \|\brho\|_1
\Big[1-\cI(\cY)\Big] \le \|\brho_\star\|_1 \Big[1 + \cI(\cY)\Big]
- r \|\brho_\star^{(o)}\|_1 .
\]
The first inequality is because $\brho_\star$ is the $\ell_1$
minimizer in (\ref{eq:F7}). Statement (\ref{eq:FG6}) follows from this
equation. $\Box$

\textbf{Proof of Corollary \ref{cor.thm4}:} We start with equation
(\ref{eq:FGPf1}) and take inner product with vector
\begin{equation}
\bu = \sum_{q \in \cS} \sigma_q \bg_q, \qquad \sigma_q = \mbox{sign}
( \rho_q - \bar \rho_{\star q}).
\label{eq:FGPf7}
\end{equation}
We obtain
\begin{equation}
  \cT_L := \Big| \sum_{q \in \cS} \rho_q \lin \bg_q,\bu \rin - \sum_{q
    \in \cS} \sum_{j \in \mathscr{S}_q} \rho_{\star j}^{(i)} \lin
  \bg_j,\bu \rin \Big| = \Big| \sum_{j \in \mathscr{S}^c} \rho_{\star
    j}^{(o)} \lin \bg_j, \bu \rin \Big| =: \cT_R,
  \label{eq:FGPf8}
\end{equation}
and proceed as in the previous proof by bounding both sides of this
equation.

For $\cT_L$ we have
\begin{align*}
  \cT_L = \Big| \sum_{q \in \cS} \sigma_q (\rho_q -\bar \rho_{\star
      q}) + \sum_{q \in \cS} \rho_q \sum_{j \in \cS, j \ne q}
  \sigma_j \lin \bg_q, \bg_j \rin - \sum_{q \in \cS} \sum_{j \in
    \mathscr{S}_q} \rho_{\star j}^{(i)} \sum_{l \in \cS, l \ne q}
  \sigma_l \lin \bg_j, \bg_l \rin \Big|,
\end{align*}
where we used definition (\ref{eq:FG8}) of the components of $\bar
\brho_\star$. We bound it as 
\begin{align}
  \cT_L &\ge \sum_{q \in \cS} |\rho_q - \bar \rho_{\star q}| - \sum_{q
    \in \cS} |\rho_q| \sum_{j \in \cS, j \ne q} |\lin \bg_q, \bg_j \rin |
  - \sum_{q \in \cS} \sum_{j \in \mathscr{S}_q} |\rho_{\star j}^{(i)}|
  \sum_{l \in \cS, l \ne q} | \lin \bg_j, \bg_l \rin | \nonumber \\ &
  \ge \|\brho - \bar \brho_{\star}\|_1 - \Big( \|\brho\|_1 +
  \|\brho_\star^{(i)} \|_1 \Big) \cI(\cY), \label{eq:FGPf9}
\end{align}
using the triangle inequality and the definition of $\sigma_q$ and
$\cI(\cY)$.

For $\cT_R$ we have by definition (\ref{eq:FGPf7}) of $\bu$ and the
triangle inequality that
\begin{align*}
  \cT_R \le \sum_{j \in \mathscr{S}^c} |\rho_{\star j}^{(o)}| \,|\lin
  \bg_j, \bg_{\cN(\vz_j)} \rin| + \sum_{j \in \mathscr{S}^c}
  |\rho_{\star j}^{(o)}| \sum_{q \in \cS, \vz_q \ne \cN(\vz_j)}
  | \lin \bg_j, \bg_q \rin |.
\end{align*}
The right hand side in this equation can be bounded as in the proof of
Theorem \ref{thm.4}, and the result is
\begin{equation}
  \cT_R \le \|\brho_\star^{(o)}\|_1 \Big[ 1-r + \cI(\cY_1) \Big].
  \label{eq:FGPf10}
\end{equation}
Now equations (\ref{eq:FGPf8})-(\ref{eq:FGPf10}) give
\begin{align*}
\|\brho - \bar\brho_\star\|_1 & \le \Big( \|\brho\|_1 +
\|\brho_\star^{(i)} \|_1 \Big) \cI(\cY) + \|\brho_\star^{(o)}\|_1
\Big[ 1-r + \cI(\cY_1) \Big] \\ & = \Big( \|\brho\|_1 + \|\brho_\star
\|_1 \Big) \cI(\cY) + \|\brho_\star^{(o)}\|_1 (1-r) \\ &\le \Big(
\|\brho\|_1 + \|\brho_\star \|_1 \Big) \cI(\cY) + \|\brho_\star\|_1
\frac{2 (1-r) \cI(\cY)}{r},
\end{align*}
with the last inequality due to Theorem \ref{thm.4}. Corollary
\ref{cor.thm4} follows from this inequality and $\|\brho_\star\|_1 \le
\|\brho\|_1$. $\Box$
\subsubsection{$\ell_1$ optimal reconstructions of clusters of sources}
\label{sect:FGProof3}
The proof of Theorem \ref{thm.6} is a slight modification of that in
section \ref{sect:FGProof1}. We begin by defining the index map
$J_\ep:\cS \to \cS_\ep$ that takes any $j \in \cS$ to $J_\ep(j)$, the
index of the point in $\cS_\ep$ at the center of the ball containing
$\vz_j$ i.e., $\vz_j \in \cB_\ep(\vz_{J_\ep(j)})$. Obviously, the
restriction of $J_\ep$ on $\cS \cap \cS_\ep$ is the identity map.

Using the definition of $\brho_\star$ and its decomposition in
$\brho_\star^{(i)}$ and $\brho_\star^{(o)}$ we obtain the equivalent
of equation (\ref{eq:FGPf1})
\begin{align}
  \sum_{q \in \cS} \rho_q\bg_q = \sum_{q\in \cS_\ep} \sum_{j \in
    \mathscr{S}_q} \rho_{\star j}^{(i)} \bg_j + \sum_{j \in
    \mathscr{S}_\ep^c} \rho_{\star j}^{(o)} \bg_j,
  \label{eq:CLP1}
\end{align}
where $\mathscr{S}_\ep^c = \{1, \ldots, N \} \setminus \displaystyle
\bigcup_{q \in \cS_\ep}\mathscr{S}_q$. We take the inner product of both
sides of this equation with vector
\begin{equation}
  \bu = \sum_{q \in \cS_\ep} \sigma_q \bg_q, \qquad \sigma_q =
  \mbox{sign} ( \bar \rho_q),
  \label{eq:CLP2}
\end{equation}
and get
\begin{equation}
  \cT_L := \Big| \sum_{q \in \cS} \rho_q \lin \bg_q,\bu \rin \Big| =
  \Big| \sum_{q\in \cS_\ep} \sum_{j \in \mathscr{S}_q} \rho_{\star
    j}^{(i)} \lin \bg_j,\bu \rin + \sum_{j \in \mathscr{S}_\ep^c}
  \rho_{\star j}^{(o)} \lin \bg_j,\bu \rin \Big| =: \cT_R,
  \label{eq:CLP3}
\end{equation}
where $\cT_L$ and $\cT_R$ denote the left and right hand side of the
equation, as before.

For $\cT_L$ we have
\begin{align*}
  \cT_L &= \Big| \sum_{q \in \cS_\ep} \sigma_q \sum_{j \in \cS \cap
    \cB_\ep(\vz_q)} \rho_j \lin \bg_j,\bg_q \rin + \sum_{q \in \cS_\ep}
  \sigma_q \sum_{j \in \cS, \vz_j \notin \cB_\ep(\vz_q)} \rho_j \lin
  \bg_j,\bg_q \rin \Big| \\ &= \Big| \sum_{q \in \cS_\ep} |\bar \rho_q|
  + \sum_{q \in \cS_\ep} \sigma_q \sum_{j \in \cS, \vz_j \notin
    \cB_\ep(\vz_q)} \rho_j \lin \bg_j,\bg_q \rin \Big| \\ & \ge
  \sum_{q \in \cS_\ep} |\bar \rho_q| - \sum_{j \in \cS} |\rho_j| \sum_{q
    \in \cS_\ep, q \ne J_\ep(j)} |\lin \bg_j,\bg_q \rin| \\ & \ge
  \sum_{q \in \cS_\ep} |\bar \rho_q| - \sum_{j \in \cS} |\rho_j|
  \cI(\cY_\ep).
\end{align*}
The equality in the second row is by definition (\ref{eq:CL2}) of the
effective source vector $\bar \rho$ and definition (\ref{eq:CLP1}) of
$\sigma_q$. The bound in the third row is by the triangle inequality
and in the last row by definition of $\cI(\cY_\ep)$. Thus, the left hand
side of (\ref{eq:CLP3}) satisfies
\begin{equation}
  \cT_L \ge \|\bar \brho\|_1 - \|\brho\|_1 \cI(\cY_\ep).
  \label{eq:CLP4}
\end{equation}

For the right hand side $\cT_R$ we have from the definition of $\bu$
and the triangle inequality
\begin{align*}
  \cT_R &= \Big| \sum_{q \in \cS_\ep} \sum_{j \in \mathscr{S}_q}
  \rho_{\star j}^{(i)} \Big[\sigma_q \lin \bg_j,\bg_q \rin + \sum_{l
      \in \cS_\ep, l \ne q}\, \sigma_l \lin \bg_j,\bg_l \rin \Big]+
  \sum_{j \in \mathscr{S}_\ep^c} \rho_{\star j}^{(o)} \sum_{l\in
    \cS_\ep}\sigma_l \lin \bg_j, \bg_l \rin \Big| \\ & \le \sum_{q \in
    \cS_\ep} \sum_{j \in \mathscr{S}_q} | \rho_{\star j}^{(i)}|
  \Big[|\lin \bg_j,\bg_q \rin| + \sum_{l \in \cS_\ep, l \ne q} |\lin
    \bg_j,\bg_q \rin|\Big] + \sum_{j \in \mathscr{S}_\ep^c}
  |\rho_{\star j}^{(o)}|\sum_{l\in \cS_\ep}|\lin \bg_j, \bg_l \rin|.
\end{align*}
In the first term we can only say that $ |\lin \bg_j, \bg_q \rin | \le
1$, because the points indexed by $\mathscr{S}_q$ are all clustered
around $\vz_q$. The sum in the second term is bounded by the
interaction coefficient of the set $\cY_\ep$, and for the last term we
have
\begin{align*}
\sum_{j \in \mathscr{S}_\ep^c} |\rho_{\star j}^{(o)}| \sum_{l\in
  \cS_\ep}| \lin \bg_j, \bg_l \rin | &= \sum_{j \in \mathscr{S}_\ep^c}
|\rho_{\star j}^{(o)}| \Big[ |\lin \bg_j, \bg_{J_\ep(j)} \rin | +
  \sum_{l \in \cS_\ep, l \ne J_\ep(j)} |\lin \bg_j,\bg_l \rin| \Big]
\\ &\le \sum_{j \in \mathscr{S}_\ep^c} |\rho_{\star j}^{(o)}| \Big[
  1-r + \cI (\cY_\ep) \Big].
\end{align*}
The upper bound on $\cT_R$ becomes
\begin{align}
  \cT_R &\le \|\brho_\star^{(i)}\|_1 \Big[1 + \cI(\cY_\ep)\Big] +
  \|\brho_\star^{(o)}\|_1 \Big[ 1 - r + \cI(\cY_\ep)\Big] \nonumber
  = \| \brho_\star \|_1 \Big[1 + \cI(\cY_\ep)\Big] - r
  \|\brho_\star^{(o)}\|_1,
    \label{eq:CLP5}
\end{align}
and substituting it in (\ref{eq:CLP3}) and using (\ref{eq:CLP4}), we
get after some rearrangement
\begin{equation}
  \|\bar \brho \|_1 - \|\brho\|_1 + |\brho\|_1
  \Big[1-\cI(\cY_\ep)\Big] \le \| \brho_\star \|_1 \Big[1 +
    \cI(\cY_\ep)\Big] - r \|\brho_\star^{(o)}\|_1.
\end{equation}
The statement of Theorem \ref{thm.6} follows from this, the assumption
that $\cI(\cY_\ep) \le 1$ and $\|\brho_\star\|_1 \le
\|\brho\|_1$. $\Box$
\subsubsection{$\ell_1$ penalty reconstructions of well separated sources}
\label{sect:FGProof2}
Before giving the proof of Theorem \ref{thm.5}, let us introduce some
notation. The sources are at points in $\cY = \{\vy_q, ~ ~ q = 1,
\ldots, s\}$ which may be off-grid, and we let $\sS$ be the set of
indexes of the grid points in the $r-$vicinity of the sources, so that
\begin{equation}
  \vz_j \in \bigcup_{q = 1}^s \cB_r (\vy_q), \qquad \forall \, j \in
  \sS.
  \label{eq:LP1}
\end{equation}
The complement of the set $\sS$ is $\sS^c = \{1,\ldots, N\} \setminus
\sS$. For any vector $\bu \in \complex^N$, we denote by $\bu_{_\sS}$
its restriction to the set $\sS$. This is a vector of length $|\sS| <
N$. We also let $\bG_{_\sS}$ be the $M \times |\sS|$ matrix with
$|\sS|$ columns $\bg_j$, for $j \in \sS$, and denote by $\cP_{_\sS}$
the orthogonal projection on the range of $\bG_{_\sS}$. It is given by
\begin{equation}
  \cP_{_\sS} = \bG_{_\sS} \bG_{_\sS}^\dagger,
  \label{eq:LP1p}
\end{equation}
where $\dagger$ denotes the pseudo-inverse.

The proof of Theorem \ref{thm.5} is based on the next two lemmas. The
first uses results in convex analysis, specifically the sub-gradient
of a convex function, defined in \cite{convex}. We need here the
sub-gradient of the $\ell_1$ norm function evaluated at $\bu \in
\complex^N$, which is shown in \cite{tropp2006just} to be any vector
in the set
\begin{equation}
  \partial \|\bu \|_1 = \{ \bxi \in \complex^N ~ ~ \mbox{s.t.}~~ \xi_i
  = \mbox{sign}(u_i) ~\mbox{if}~ u_i \ne 0 ~~\mbox{and} ~~ |\xi_i| \le
  1 ~\mbox{if}~ u_i = 0 \}.
  \label{eq:LP2}
\end{equation}
The second lemma estimates the Lagrange multiplier $\gamma$ needed to
prove  the theorem.
\begin{lemma}
  \label{lem1.thm5}
  Let $\brho_\star$ minimize the augmented Lagrangian $\sL(\brho)$ defined
  in (\ref{eq:F7b}),
  over vectors supported in $\sS$, with $\brho_{\star_\sS}$ its
  restriction to $\sS$. Then, there exists a sub-gradient vector $\bxi
  \in \partial \|\brho_{\star_\sS}\|_1$ such that
  \begin{equation}
    \bG_{_\sS}^{\tiny H} \bG_{_\sS} \big(\brho_{\star_\sS} -
    \bG_{_\sS}^\dagger\bd \big) + \gamma \bxi = {\bf 0},
    \label{eq:LP4}
  \end{equation}
  where the index $H$ denotes the Hermitian adjoint of $\bG_{_\sS}$, a
  matrix in $\complex^{|\sS| \times N}$. Moreover, if we let $\cS
  \subset \sS$ be the set of $s$ grid points that are nearest the
  locations $\vy_j$ of the $s$ sources, we get
  \begin{equation}
    \bG_{_\cS}^H \bG_{_\sS} \big(\brho_{\star_\sS} - \bG_{_\sS}^\dagger\bd
    \big) + \gamma \bxi_{_\cS} = {\bf 0}.
    \label{eq:LP5}
  \end{equation}
\end{lemma}

\textbf{Proof:} For any $\brho$ supported in $\sS$ we have
$\|\brho\|_1 = \|\brho_{_\sS} \|_1$, and using Pythagora's theorem
\begin{align*}
\|\bG \brho - \bd\|_2^2 &= \|\bG_{_\sS} \brho_{_\sS} - \cP_{_\sS}
\bd\|_2^2 + \|\bd - \cP_{_\sS} \bd\|_2^2.
\end{align*}
Therefore
$
\mathscr{L}(\brho) = \mathscr{L}_{_\sS}(\brho_{_\sS}) + \frac{1}{2}\|\bd -
\cP_{_\sS} \bd\|_2^2,
$
with $\sL_{_\sS}$ defined on vectors of length $|\sS|$, 
\begin{equation}
\mathscr{L}_{_\sS}(\bx) = \frac{1}{2} \|\bG_{_\sS} \bx - \cP_{_\sS}
\bd \|_2^2 + \gamma \|\bx \|_1, \qquad \forall \, \bx \in
\complex^{|\sS|}.
\end{equation}
We conclude that $\brho_{{\star_\sS}}$ is the minimizer of
$\mathscr{L}_{_\sS}$. Then, results in convex analysis
\cite{convex,tropp2006just} imply that ${\bf 0}$ must be an element of
the sub-gradient of $\mathscr{L}_{_\sS}$. Equivalently, there exists a
vector $\bxi \in \partial \|\brho_{\star_{\sS}}\|_1$ satisfying
(\ref{eq:LP4}), where we use the expression of the projection
$\cP_{_\sS}$. Equation (\ref{eq:LP5}) is just the restriction of
equation (\ref{eq:LP4}) to the rows indexed by $\cS$. $\Box$

\vspace{0.05in}
\begin{lemma}
  \label{lem2.thm5}
  Let $\brho_\star$ be the minimizer of $\mathscr{L}(\brho)$ over
  vectors supported in $\sS$. If $\gamma$ satisfies
  \begin{equation}
    \sqrt{2 r} \|\brho_{\star_\sS} - \bG_{_\sS}^\dagger \bd\|_1 +
    \|\bG_{_{\sS^c}}^H (\cP_{_\sS} \bd - \bd) \|_{\infty} < \gamma
    \Big[ 1 - \max_{j \in \sS^c} \Big| \lin \bxi_{_\cS},
      \bG_{_\cS}^\dagger \bg_j \rin \Big| \Big],
    \label{eq:LP6}
  \end{equation}
  with $\bxi$ as in Lemma \ref{lem1.thm5}, $\brho_\star$ is the global
  minimizer of $\mathscr{L}(\brho)$ over $\complex^N$.
\end{lemma}

\textbf{Proof:} To prove the lemma we show that any perturbation of
$\brho_\star$ by a vector that is not supported in $\sS$ leads to an
increase of the objective function $\sL$. This implies that
$\brho_\star$ is a local minimizer of $\sL$ in $\complex^N$. That
$\brho_\star$ is the global minimizer follows from the convexity of
$\mathscr{L}$.

Consider an arbitrary vector $\bv \in \complex^N$ and decompose it as
\begin{equation}
  \label{eq:LP00}
\bv = \bu + \bw, \qquad \rm{supp} \, \bu \subset \sS, ~ ~ \rm{supp} \,
\bw \subset \sS^c.
\end{equation}
For small and positive $\ep$ we have from definition (\ref{eq:F7b})
and the disjoint support of $\bu$ and $\bw$ that
\begin{align}
  \sL(\brho_\star + \ep \bv) - \sL(\brho_\star + \ep \bu) = \ep \Big[
    \mbox{real} \Big( \lin \bG \brho_\star - \bd, \bG \bw \rin \Big) +
    \gamma \|\bw\|_1 \Big]+  \ep^2 \Big[\frac{1}{2} \|\bG
    \bw\|_2^2 + \mbox{real} \Big( \lin \bG \bu, \bG \bw \rin
    \Big)\Big],
  \label{eq:LP7a}
\end{align}
with the first term dominating the second for $\ep \ll 1$. We write it
in terms of the components $w_j$ of $\bw$ as
\begin{align}
  \lin \bG \brho_\star - \bd, \bG \bw \rin &= \sum_{j \in \sS^c} \lin
  \bG_{_\sS} \brho_{\star_\sS} - \bd, \bg_j\rin w_i =
  \sum_{j \in \sS^c} \lin \bG_{_\sS} \brho_{\star_\sS} - \cP_{_\sS}
  \bd, \bg_j\rin w_j - \sum_{j \in \sS^c} \lin \bd - \cP_{_\sS} \bd,
  \bg_j\rin w_j,
  \label{eq:LP7}
\end{align}
where we used that $\brho_\star$ is supported in $\sS$. We estimate
next the two sums in the right hand side.

For the terms in the first sum we have
\begin{align}
  \hspace{-0.1in} \lin \bG_{_\sS} \brho_{\star_\sS} - \cP_{_\sS} \bd, \bg_j\rin
  =\lin \bG_{_\sS} \brho_{\star_\sS} - \cP_{_\sS} \bd, \cP_{_\cS}
  \bg_j\rin + \lin \bG_{_\sS} \brho_{\star_\sS} - \cP_{_\sS} \bd,
  (I-\cP_{_\cS}) \bg_j\rin,
  \label{eq:LP9}
\end{align}
where $\cP_{_\cS} = \bG_{_\cS} \bG_{_\cS}^\dagger$ is the orthogonal
projection on the range of the full rank \footnote{That $\bG_{_\cS}$
  is full rank follows from the assumption $\cI(\cY) < 1/2$. Matrix
  $\bG_{_\sS}$ is not full rank because its columns are associated to
  nearby points in the vicinity of the sources, as stated in
  (\ref{eq:LP1}).} matrix $\bG_{_\cS}$, with pseudo-inverse
\begin{equation}
    \bG_{_\cS}^\dagger =\Big(
   \bG_{_\cS}^H \bG_{_\cS} \Big)^{-1} \bG_{_\cS}^H.
   \label{eq:LP8}
\end{equation}
Substituting the expression of $\cP_{_\cS}$ and $\cP_{_\sS}$ in
(\ref{eq:LP9}) we get
\begin{align}
  \lin \bG_{_\sS} \brho_{\star_\sS} - \cP_{_\sS} \bd, \bg_j\rin &=
  \lin \bG_{_\cS}^H \bG_{_\sS} \big(\brho_{\star_\sS} - \bG_{_\sS}^\dagger \bd
  \big), \bG_{_\cS}^\dagger \bg_j \rin + \lin \bG_{_\sS}
  \Big(\brho_{\star_\sS} - \bG_{_\sS}^\dagger \bd\Big), (I-\cP_{_\cS})
  \bg_j\rin \nonumber \\ &= - \gamma \lin
  \bxi_{_\cS},\bG_{_\cS}^\dagger \bg_j \rin + \lin \bG_{_\sS} \balpha,
  (I-\cP_{_\cS}) \bg_j\rin,
  \label{eq:LP10}
\end{align}
with the second inequality following from Lemma \ref{lem1.thm5}, and
notation $ \balpha = \brho_{\star_\sS} - \bG_{_\sS}^\dagger \bd.  $
Let us write explicitly the second term in (\ref{eq:LP10})
\begin{align}
   \lin \bG_{_\sS} \balpha, (I-\cP_{_\cS}) \bg_j\rin = \sum_{l\in \sS}
   \lin \bg_l, (I-\cP_{_\cS}) \bg_j\rin \alpha_l, \qquad \rm{for} ~
   j \in \sS^c.
   \label{eq:PL11}
\end{align}
Since $l \in \sS$, we have by definition (\ref{eq:LP1}) that there
exists $q \in \cS$ such that $\vz_l \in \cB_r(\vy_q)$, and we can
decompose $\bg_l$ in two parts: $\bg_l^\parallel$ which is along
$\bg_q$ and $\bg_l^\perp$ which is orthogonal to it,
\[
\bg_l = \bg_l^{\parallel} + \bg_l^\perp, \qquad
\bg_l^{\parallel} = \lin \bg_q, \bg_l \rin \bg_q, \quad
\bg_l^\perp = \bg_l - \bg_l^{\parallel}.
\]
Clearly $\bg_l^\parallel$ is in the range of $\bG_{_\cS}$, so it is
orthogonal to $(I-\cP_{_\cS})\bg_j$, and the terms in equation
(\ref{eq:PL11}) satisfy 
\begin{align}
   \big|\lin \bg_l, (I-\cP_{_\cS}) \bg_j\rin\big| &= \big| \lin
   \bg_l^\perp, (I-\cP_{_\cS}) \bg_j\rin \big| \le \|\bg_l^\perp\|_2
   \|(I-\cP_{_\cS}) \bg_j\|_2 \le \|\bg_l^\perp\|_2.
   \label{eq:PL12}
\end{align}
Moreover, by Pythagora's theorem 
\begin{equation}
\|\bg_l^\perp\|_2^2 = 1 - \|\bg_l^\parallel\|_2^2 = 1 - \big| \lin
\bg_q, \bg_l \rin \big|^2 < 1 - (1-r)^2 = 2r-r^2 < 2r,
\label{eq:PL13}
\end{equation}
where the first inequality follows $\vz_l \in \cB_r(\vy_q)$ i.e.,
$\D(\vz_l,\vz_q) = 1 - \big|\lin \bg_q,\bg_l \rin \big| < r.$
Gathering the results (\ref{eq:LP10})-(\ref{eq:PL13}) and using the
triangle inequality we get the following bound on the first
sum in (\ref{eq:LP7})
\begin{equation}
  \big| \sum_{j \in \cS^c} \lin \bG_{_\sS} \brho_{\star_\sS} -
  \cP_{_\sS} \bd, \bg_j \rin w_j \big| \le \Big[\gamma \max_{j \in
      \sS^c} \big| \lin \bxi_{_\cS}, \bG_{_\cS}^\dagger \bg_j \rin
    \big| + \sqrt{2 r} \|\brho_{\star_\sS}-\bG_{_\sS}^\dagger \bd
    \|_1\Big] \|\bw\|_1.
  \label{eq:PL14}
\end{equation}

For the second sum in (\ref{eq:LP7}) we have
\begin{align}
  \big| \sum_{j \in \sS^c} \lin \bd - \cP_{_\sS} \bd, \bg_j \rin w_j
  \big| &\le \max_{l \in \sS^c} \big| \lin \bd - \cP_{_\sS} \bd, \bg_l
  \rin \big| \sum_{j \in \sS^c} |w_j| 
  = \|\bG_{_{\sS^c}}^H (\bd - \cP_{_\sS} \bd) \|_{\infty} \|\bw\|_1,
  \label{eq:PL15}
\end{align}
and putting together the results (\ref{eq:PL14})-(\ref{eq:PL15}) we
obtain from (\ref{eq:LP7}) that
\begin{align}
  \big| \lin \bG \brho_\star - \bd, \bG \bw \rin \big| \le \Big[
    \gamma \max_{j \in \sS^c} \big| \lin \bxi_{_\cS},
    \bG_{_\cS}^\dagger \bg_j \rin \big| + \sqrt{2 r}
    \|\brho_{\star_\sS}-\bG_{_\sS}^\dagger \bd \|_1 +
    \|\bG_{_{\sS^c}}^H (\bd - \cP_{_\sS} \bd) \|_{\infty}
    \Big]\|\bw\|_1.
\end{align}
This estimate and the triangle inequality give that the $\ep$ term in
(\ref{eq:LP7a}) is positive for $\gamma$ satisfying
\begin{align*}
  \gamma > \gamma \max_{j \in \sS^c} \big| \lin \bxi_{_\cS},
    \bG_{_\cS}^\dagger \bg_j \rin \big| + \sqrt{2 r}
    \|\brho_{\star_\sS}-\bG_{_\sS}^\dagger \bd \|_1 + 
    \|\bG_{_{\sS^c}}^H (\bd - \cP_{_\sS} \bd) \|_{\infty},
\end{align*}
as assumed in the lemma. Then, the perturbation of $\brho_\star$ by
the arbitrary vector (\ref{eq:LP00}) increases the objective function
$\sL$ and the lemma follows. $\Box$

\vspace{0.05in}
To complete the proof of the theorem it remains to show that we can
find a positive $\gamma$ as in Lemma \ref{lem2.thm5}. We begin with
the estimate
\begin{align}
  \max_{j \in \sS^c} \big| \lin \bxi_{_\cS}, \bG_{_\cS}^\dagger \bg_j
  \rin \big| &= \max_{j \in \sS^c} \left| \lin \bxi_{_\cS},
  \big(\bG_{_\cS}^H \bG_{_\cS}\big)^{-1} \bG_{_\cS}^H \bg_j \rin
  \right| \nonumber \\ &\le \max_{j \in \sS^c} \big\| \bxi_{_\cS}
  \|_{\infty} \big \| \Big(\bG_{_\cS}^H \bG_{_\cS}\Big)^{-1}
   \bG_{_\cS}^H \bg_j \big \|_1 \nonumber \\ &\le
  \max_{j \in \sS^c}\big \| \big(\bG_{_\cS}^H \bG_{_\cS}\big)^{-1}
  \big\|_{1,1} \big\| \bG_{_\cS}^H \bg_j \big \|_1,
  \label{eq:PL16}
\end{align}
where we used (\ref{eq:LP8}) and definition (\ref{eq:LP2}) of the
sub-gradient of the $\ell_1$ norm. Now since $j \in \sS^c$, $\vz_j$ is
outside every ball of radius $r$ centered at a source or,
equivalently, $ |\lin \bg_q, \bg_j \rin| \le 1 -r$, for all $q \in
\cS.$ This and the definition of $\cI(\cY)$ give
\begin{equation}
  \big\| \bG_{_\cS}^H \bg_j \big \|_1 = \sum_{q \in \cS} \big| \lin
  \bg_q, \bg_j \rin \big| \le 1 - r + \cI(\cY), \qquad \forall \, j
  \in \sS^c.
  \label{eq:PL17}
\end{equation}
Moreover, for any vector $\bu$ supported on $\cS$ we have
\begin{align*}
  \Big\| \bG_{_\cS}^H \bG_{_\cS} \bu \big\|_1 &= \sum_{q \in \cS}
  \big| \sum_{j \in \cS} \lin \bg_q,\bg_j \rin u_j \big| \\ &= \sum_{q
    \in \cS} \big| u_q + \sum_{j \in \cS, j \ne q} \lin \bg_q,\bg_j
  \rin u_j \big| \\ &\ge \sum_{q \in \cS} \big| u_q | \Big[ 1 -
    \sum_{j \in \cS, j \ne q} \big|\lin \bg_q,\bg_j \rin \big| \Big] \\
  & \ge \Big[1 - \cI(\cY)\Big] \|\bu\|_1,
\end{align*}
so the operator norm of the inverse of $\bG_{_\cS}^H \bG_{_\cS}$,
which exists because $\bG_{_\cS}$ is full rank, satisfies
\begin{equation}
  \Big\| \big( \bG_{_\cS}^H \bG_{_\cS}\big)^{-1} \Big\|_{1,1} \le
  \Big[1 - \cI(\cY)\Big]^{-1}.
  \label{eq:PL18}
\end{equation}
Putting together (\ref{eq:PL16})-(\ref{eq:PL18}) we obtain that
\begin{equation}
\max_{j \in \sS^c} \big| \lin \bxi_{_\cS}, \bG_{_\cS}^\dagger \bg_j
\rin \big| \le \frac{1- r + \cI(\cY)}{1 - \cI(\cY)} < 1,
\label{eq:PL19}
\end{equation}
where the second inequality is by the assumption of the theorem that
$r > 2 \cI(\cY)$. This shows that the right hand side in equation
(\ref{eq:LP6}) in Lemma \ref{lem2.thm5} is positive.

Finally, we show that the left hand side in equation (\ref{eq:LP6}) is
bounded independent of $\gamma$. Clearly, the term $\|\bG_{_{\sS^c}}^H
(\cP_{_\sS} \bd - \bd) \|_{\infty}$ does not depend on $\gamma$. It is
due to the modeling error that is small when the grid is fine enough
and the additive noise, that may cause $\bd$ to lie outside the range
of $\bG_{_\sS}$ i.e., $\bd \ne \cP_{_\sS} \bd.$ To bound the first
term in (\ref{eq:LP6}) note that $ \bu = \bG_{_\sS}^\dagger \bd $ is
the minimum $\ell_2$ norm solution of $\bG_{_\sS} \bu = \bd$ if it
exists, or otherwise the minimizer of the least squares misfit
$\|\bG_{_\sS} \bu - \bd\|_2$. Then, since $\brho_{\star_\sS}$
minimizes $\sL_{_\sS}$, we have
\begin{align*}
  \sL_{_\sS}(\brho_{\star_\sS}) = \frac{1}{2} \|\bG_{_\sS}
  \brho_{\star_\sS} - \bd\|_2^2 + \gamma \|\brho_{\star_\sS} \|_1 \le
 \sL_{_\sS}(\bu) = \frac{1}{2} \|\bG_{_\sS}
 \bu - \bd\|_2^2 + \gamma \|\bu \|_1,
\end{align*}
and using that $\|\bG_{_\sS} \bu - \bd\|_2 \le \| \bG_{_\sS}
\brho_{\star_\sS} - \bd \|_2$, we get $ \|\brho_{\star_\sS} \|_1 \le
\| \bu\|_1.  $ Consequently,
\begin{equation}
  \|\brho_{\star_\sS} - \bG_{_\sS}^\dagger \bd \|_1 =
  \|\brho_{\star_\sS} - \bu \|_1 \le \|\brho_{\star_\sS}\|_1 +
  \|\bu\|_1 = 2 \|\bu_1\|,
\end{equation}
independent of $\gamma$. We conclude that we can find $\gamma > 0$ as
in Lemma \ref{lem2.thm5}, and therefore complete the proof of the
theorem. $\Box$

\section{Summary}
\label{sect:sum}
We presented a resolution study of sensor array imaging of a sparse
scene of point sources or scatterers. The setup is in the paraxial
regime, where the array aperture is small with respect to the distance
to the imaging region. The imaging is done with two sparsity promoting
optimization methods: $\ell_1$ optimization (basis pursuit) and
$\ell_1-$penalty. The latter deals with noise and modeling errors.
Our resolution analysis takes into account the sparse support of the
unknowns. In case that they lie on the imaging grid, we obtained
conditions on the grid size that guarantee their exact and unique
recovery for noiseless data.  This is for both single frequency and
broad-band regimes, and the results show the benefit of having
multiple frequency measurements. In case that the unknowns lie
off-grid, we studied imaging on fine grids that mitigate the modeling
error. We showed that when the unknowns are located at sufficiently
far apart points in the scene, or they lie in well separated clusters,
the results of imaging with sparsity promoting optimization are
useful. The support of the reconstruction is near that of the unknowns
and its locally averaged amplitudes approximates the true ones.

\section*{Acknowledgments}
This work was partially supported by AFOSR Grant FA9550-15-1-0118. LB
also acknowledges support from the ONR Grant N00014-14-1-0077.

 \appendix

\section{Numerical setup}
\label{ap:numeric}
The simulations are for an aperture $a = 25 \la$, range $L = 1000
\la$. We used different sizes of imaging grids such as $W_N = 10 \times 10 \times 20$ or $W_N = 5 \times 64 \times 64$, for a given
mesh size $\vbh = (h,h,h_3)$.

The $\ell_1$ optimization is solved with the package \cite{cvx}. For
noisy data we solve (\ref{eq:F7b}) with $\gamma$ chosen to be close to
the noise level. We find the results to be very similar to those
obtained from the constrained optimization
\[
\min_{\brho \in \complex^N} \|\brho\|_1 \quad \mbox{s.t.} \quad
\|\bG \brho-\bd\|_2 \le ~ \mbox{noise level}.
\]
Because the simulations give only an approximation $\widetilde
\brho_\star$ of the minimizer $\brho_\star$, we threshold the results
at $1\%$ of the maximum entry in absolute value. We say that $\brho$
is recovered numerically if $ {\|\widetilde \brho_\star -
  \brho\|_{\infty}}/{ \|\brho\|_{\infty}} < 1\%.  $

\section{Derivation of the paraxial model}
\label{ap:parax_deriv}
We have for $\vx = (\bx,0)$ in the array and $\vz = (\bz,z_3)$ in the
imaging region that 
\begin{equation}
|\vx-\vy| = L \Big[ 1 + O \Big(\frac{D_3}{L}\Big) + O
  \Big(\frac{a^2}{L^2}\Big) \Big]
\label{eq:C0}
\end{equation}
so we can approximate the geometrical spreading factor in the Green's function by 
\begin{equation}
\frac{1}{4 \pi |\vx-\vy|} \approx \frac{1}{4 \pi L}.
\label{eq:C01}
\end{equation}
The phase is given by
\begin{equation}
k |\vx-\vy| = k z_3 + \frac{k|\bx-\bz|^2}{2 z_3} + \cE(\vx,\vy),
\label{eq:C1}
\end{equation}
with remainder 
\begin{equation*}
\cE(\vx,\vy) = - \frac{k |\bx-\bz|^4}{8 z_3^3} +
O\Big(\frac{k|\bx-\bz|^6}{z_3^5}\Big).
\end{equation*}
The scaling assumptions (\ref{eq:P4}) make most of the
terms in $\cE$ negligible, except for those that are independent of
$\vz$, which cancel in the product of the Green's functions in the
left hand side of equation (\ref{eq:P8}). Thus we write 
\begin{equation}
\cE(\vx,\vy) \approx -\frac{k |\bx|^4}{8 L^3} +
O\Big(\frac{k|\bx|^6}{L^5}\Big).
\label{eq:C2}
\end{equation}
We also have
\begin{align*}
\frac{k|\bx-\bz|^2}{2 z_3} = \frac{k |\bx|^2}{2 z_3} - \frac{k \bx
  \cdot \bz}{z_3} + \frac{k |\bz|^2}{2 z_3},
\end{align*}
with the last term negligible by assumption (\ref{eq:P4}). Simplifying
further, we get
\begin{align}
\frac{k|\bx-\bz|^2}{2 z_3} \approx \frac{k |\bx|^2}{2 L} - \frac{k
  |\bx|^2 z_3}{2 L^2} - \frac{k \bx \cdot \bz}{L} +
O\Big(\frac{a^2}{\la L} \frac{D_3^2}{L^2} \Big) + O \Big( \frac{a D
  D_3}{\la L^2}\Big),
\label{eq:C3}
\end{align}
with the remainder negligible by (\ref{eq:P4}).

Now let $\vx = \vx_r$ and obtain from (\ref{eq:C01})-(\ref{eq:C3})
that
\begin{align*}
\hat G(\om,\vx_r,\vz_j) \overline{\hat G(\om,\vx_r,\vz_q)} \approx
\frac{e^{i k (z_{3,j}-z_{3,q})}}{(4 \pi L)^2} e^{- \frac{i k
    |\bx_r|^2(z_{3,j}-z_{3,q})}{2 L^2} - \frac{i k \bx_r \cdot
    (\bz_j-\bz_q)}{L}}.
\end{align*}
Equation (\ref{eq:P8}) follows by summing over $r$. $\Box$

\section{Derivation of the broad-band paraxial model}
\label{ap:BBparax_deriv}
Using the parabolic scaling, we find as in appendix
\ref{ap:parax_deriv} that
\begin{equation}
 \hat G(\om_j,\vx_r,\vz_q) \overline{\hat G(\om_j,\vx_r,\vz_l)} \approx
\frac{e^{i k_j (z_{3,q}-z_{3,l})}}{(4 \pi L)^2} e^{- \frac{i k_j
    |\bx_r|^2(z_{3,q}-z_{3,l})}{2 L^2} - \frac{i k_j \bx_r \cdot
    (\bz_q-\bz_l)}{L}},
\label{eq:BBP1}
\end{equation}
where $k_j = \om_j/c$. Moreover, the phase $\Phi$ in the right hand
side of (\ref{eq:BBP1}) is
\begin{align*}
  \Phi = k_j (z_{3,q} - z_{3,l}) - \frac{k_o
    |\bx_r|^2(z_{3,q}-z_{3,l})}{2 L^2} - \frac{k_o \bx_r \cdot
    (\bz_q-\bz_l)}{L} +  O\Big(\frac{B}{c} \frac{a D}{L}
  \Big) + O\Big(\frac{B}{c} \frac{a^2 D_3}{L^2}\Big),
\end{align*}
with negligible remainder by assumption (\ref{eq:BB2}). Result
(\ref{eq:BB3}) follows after substituting the approximation of $\Phi$
in (\ref{eq:BBP1}), multiplying with the Gaussian pulse, and summing
over the frequencies and the receivers. We also remove the phase 
$k_o (z_{3,q}-z_{3,l})$ in equation (\ref{eq:BB3}). $\Box$

\section{Proof of estimate \ref{eq:PP10}}
\label{ap:contour}
Simplifying notation as
\[
r_o = \frac{\beta-\eta}{2 \sqrt{\eta}} > \frac{\alpha}{2
  \sqrt{\eta}} > 0, \qquad r_1 = \frac{\beta+ \eta}{2
  \sqrt{\eta}} = r_0 + \sqrt{\eta},
\]
and choosing a contour defined by the line segments at angle $0$ and
$\pi/4$ with the real axis, and the circular arcs at radius $r_0$ and
$r_1$, we obtain
\begin{align}
\int_{r_0}^{r_1} \hspace{-0.05in}dt \, e^{i t^2}
=-\hspace{-0.02in}\int_0^{\pi/4} \hspace{-0.08in}d \theta \left[ r_1
  e^{-r_1^2 \sin (2 \theta) + i [\theta + r_1^2 \cos(2\theta)]} - r_0
  e^{-r_0^2 \sin (2 \theta) + i [\theta + r_0^2 \cos(2\theta)]}
  \right]  + \sqrt{i} \int_{r_o}^{r_1} dr\,
e^{-r^2}. \label{eq:apC}
\end{align}
For the last term we have the estimate
\begin{align}
\int_{r_0}^{r_1} dr\, e^{-r^2} &= \int_0^{\sqrt{\eta}} ds \,
e^{-(r_o+s)^2} 
\le e^{-r_0^2} \int_0^\infty ds \, e^{-s^2}
= \frac{\sqrt{\pi}}{2} e^{-r_o^2},
\end{align}
where we changed variables as $r = r_o + s$.  Moreover,
\begin{equation}
  \label{eq:apC2}
\left| \int_0^{\pi/4} \hspace{-0.08in}d \theta \, r_0 \, e^{-r_0^2
  \sin (2 \theta) + i [\theta + r_0^2 \cos(2\theta)]} \right| \le
\int_0^{\pi/4} \hspace{-0.08in}d \theta \, r_o e^{- 4 r_o^2 \theta
  /\pi} = \frac{\pi}{4 r_0} \Big(1-e^{-r_o^2}\Big),
\end{equation}
and similar for the other integral over $\theta$, because $ \sin (2
\theta) \ge 4 \theta/\pi$, for all $\theta \in (0, \pi/4).  $ The
estimate (\ref{eq:PP10}) follows from (\ref{eq:apC})-(\ref{eq:apC2})
and the triangle inequality
\begin{align*}
\left|\int_{r_o}^{r_1} \hspace{-0.05in}dt \, e^{i t^2}\right| &\le
\frac{\sqrt{\pi}}{2} e^{-r_o^2} + \frac{\pi}{4 r_1}\Big(1-e^{-r_1^2}\Big) +
\frac{\pi}{4 r_0}\Big(1-e^{-r_0^2}\Big) \\ & \le \frac{\pi}{2 r_0} +
\frac{\sqrt{\pi}}{2}e^{-r_0^2} <  \frac{(\pi+1) \sqrt{\eta}}{\alpha}.
\end{align*}
The second inequality is because $r_1 > r_0$, the third inequality is
by the definition of $r_0$ and the last inequality is because
$ e^{-x^2/4} < 2/(\sqrt{\pi} x)$ for any $x > 0$.  $\Box$

\bibliographystyle{siam} \bibliography{SPARSE.bib}

\end{document}